\documentclass[11pt]{article}
\usepackage{amstext}
\usepackage{amsmath}
\usepackage{amsfonts}
\usepackage{amssymb}
\usepackage{amsxtra}
\usepackage{srcltx}
\usepackage[left=1.5in,right=1in,top=1in,bottom=1in]{geometry}

\usepackage{oldgerm}
\usepackage{textcomp}
\usepackage{amsfonts}
\usepackage{amsmath}
\usepackage{amssymb}
\usepackage{amsbsy}  
\usepackage{amsthm}
\usepackage[mathscr]{eucal}
\usepackage{color}

\usepackage{ifpdf}
\usepackage{epic}

\ifpdf
  \usepackage{eepicemu}
  \usepackage[pdftex]{graphicx}
\else
  \usepackage{eepic}
  \usepackage{graphicx}
\fi

\usepackage[pagebackref]{hyperref}  
\hypersetup{%
          colorlinks=true,
          linkcolor=webbrown,
          filecolor=webbrown,
          citecolor=webgreen,
          breaklinks=true}
\ifpdf
\hypersetup{pdftex,
            pdftitle={On The Optimal Stopping Problem For One--Dimensional Diffusions},
            pdfauthor={Savas Dayanik},
            pdfstartview=FitH, 
            bookmarksopen=true,
            bookmarksnumbered=true
}
\else
\hypersetup{dvips}
\fi

\definecolor{webgreen}{rgb}{0,.5,0}
\definecolor{webbrown}{rgb}{.8,0,0}
\definecolor{emphcolor}{rgb}{0.95,0.95,0.95}

\newcommand{\hreft}[1]{\hyperref[#1]{Theorem~\ref*{#1}}}
\newcommand{\hreftb}[1]{\hyperref[#1]{Table~\ref*{#1}}}
\newcommand{\hrefl}[1]{\hyperref[#1]{Lemma~\ref*{#1}}}
\newcommand{\hrefc}[1]{\hyperref[#1]{Corollary~\ref*{#1}}}
\newcommand{\hrefp}[1]{\hyperref[#1]{Proposition~\ref*{#1}}}
\newcommand{\hreff}[1]{\hyperref[#1]{Figure~\ref*{#1}}}
\newcommand{\hrefr}[1]{\hyperref[#1]{Remark~\ref*{#1}}}
\newcommand{\hrefs}[1]{\hyperref[#1]{Section~\ref*{#1}}}
\newcommand{\hrefa}[1]{\hyperref[#1]{Appendix~\ref*{#1}}}

\newcommand{\R}{\mathbb R}
\newcommand{\N}{\mathbb N}
\renewcommand{\P}{\mathbb{P}^{\vec{\pi}}}

\newcommand{\E}{\mathbb{E}^{\vec{\pi}}}

\newcommand{\Fb}{\mathbb F}
\newcommand{\F}{\mathcal F}
\newcommand{\G}{\mathcal G}

\renewcommand{\emptyset}{\varnothing}
\newcommand{\vp}{\vec{\pi}}
\newcommand{\vP}{\vec{\Pi}}
\newcommand{\vx}{\vec{x}}
\newcommand{\eps}{\varepsilon}

\newtheorem {cor}{Corollary}[section]

\newtheorem{lemma}{Lemma}[section]
\newtheorem{remark}{Remark}[section]
\newtheorem{proposition}{Proposition}[section]



\numberwithin{equation}{section}

\usepackage{natbib}

\title{\LARGE \bf
\sc Quickest Detection for a Poisson Process with a Phase-type
Change-time Distribution
\thanks{The first author is
supported in part by the National Science Foundation, under grant
DMS-0604491.}}

\author{Erhan Bayraktar
and Semih Sezer
\thanks{E. Bayraktar and S. Sezer are with the Department of Mathematics, University
of Michigan, Ann Arbor, MI 48109, USA, email: \{erhan, sezer\}
@umich.edu} }
\date{}

\begin{document}
\maketitle
\begin{abstract}
We consider a change detection problem in which the arrival rate
of a Poisson process changes suddenly at some unknown and
unobservable disorder time. It is assumed that the prior
distribution of the disorder time is known. The objective is to
detect the disorder time with an online detection rule (a stopping
time) in a way that balances the frequency of false alarm and
detection delay. So far in the study of this problem, the prior
distribution of the disorder time is taken to be exponential
distribution for analytical tractability. Here, we will take the
prior distribution to be a phase-type distribution, which is the
distribution of the absorption time of a continuous time Markov
chain with a finite state space. We find an optimal stopping rule
for this general case and give a numerical algorithm
that calculates the parameters of $\eps$-optimal strategies for
any $\eps>0$. We illustrate our findings on two examples.
\end{abstract}

\section{Introduction}

Suppose that our observations come from a Poisson process
$X=\left\{X_t: t\geq 0\right\}$ whose arrival rate changes from
$\lambda_0$ to $\lambda_1$ at some random time $\Theta$. The
\emph{disorder time} $\Theta$ is unobservable but its prior
distribution is known. We assume that the prior distribution of
$\Theta$ is a phase-type distribution. This is the distribution of
the time of death (absorption) of a non-conservative Markov process
$M=\left\{M_t: t\geq 0\right\}$, whose state space is finite and
includes a single absorbing state. Our problem is to find an alarm
time $\tau$ which depends only on the past and the present
observations and rings as soon as $\Theta$ occurs. Since $\Theta$
is unobservable a detection rule $\tau$ will make false alarms or
have detection delays. We will find a rule that optimally balances
these two. We will choose a Bayesian risk that penalizes the sum
of the frequency of false alarm and a multiple of detection delay
as in \cite{MR2003i:60071}.

So far in the literature of continuous time Bayesian quickest
detection problems, the distribution of the disorder time $\Theta$
is always taken to be exponential distribution for analytical
tractability, see e.g. \cite{galchuk71}, \cite{davis76}, \cite{Sh78}, 
\cite{MR2001m:62090,MR2003i:60071}, \cite{MR2002b:62088}, 
\cite{MR1985648}, \cite{BD04,bdk05}, \cite{ds}, \cite{BD03}. 
The disorder time in the works cited above is modeled as the first arrival time of a Poisson process that we do not observe.
We will change the assumption on the nature of the arrivals
 for broader applicability and we will
solve the Poisson disorder problem with a phase type disorder
distribution. This seems to strike a balance between generality
and tractability. Indeed, any positive distribution may be
approximated arbitrarily closely by phase-type distributions. See
\cite{neuts} for this and other properties of this class of
distributions.

Let  $\{1,\cdots,n,\Delta\}$ denote the state space of $M$ where $\Delta$ is absorbing and the rest of the states are transient. To solve the Poisson disorder problem, we first show
that it is equivalent to an optimal stopping problem for an $n+1$
dimensional piece-wise deterministic Markov process
$\vec{\Pi}_t=[\Pi^{(1)}_t, \cdots \Pi^{(n)}_t,\Pi_t]$, $t \geq 0$,
whose $i$th coordinate is the posterior probability
$\Pi^{(i)}_t=\P\left\{M_t=i\right\}$ that the Markov chain $M$ is
in state $i$ given the past observations $\F_t=\sigma\{X_s: 0 \leq
s \leq t\}$ of $X$. The process $\Pi_t$, $t \geq 0$, is the
posterior probability that the disorder has already occurred. All
of the coordinates are driven by the same point process. We show
that the optimal stopping time (of the filtration
$\mathbb{F}=\{\F_t\}_{t \geq 0}$) is the hitting time of the
process $\vP$ to some closed convex set $\Gamma$ with
non-empty interior. We describe a numerical algorithm that
approximates the optimal Bayes risk within any given positive
error margin. Among the outputs of this algorithm are boundary curves
that characterize $\eps-$optimal stopping times. Once these curves
are determined the only thing an observer has to do is to ring the
alarm as soon as $\vP$, which is completely determined by the
observations of $X$, crosses one of these boundaries, continuously or via a jump. To see the
efficacy of the numerical algorithm we use it to approximate the
minimum Bayes risk when the prior distribution of the disorder
time has Erlang 
or Hypergeometric distribution with two non-absorbing states.

The rest of the paper is organized as follows: In
Section~\ref{sec:prob-stat}, we give a precise description of the
problem and show that it is equivalent to solving an optimal
stopping problem for 
the process $\vP$. In
Section~\ref{sec:sequential}, we show that the minimum Bayes risk
can be uniformly approximated by a sequence of functions that can
be constructed via an iterative application of an integral operator to
the terminal penalty of the optimal stopping problem described in
Section~\ref{sec:prob-stat}. A similar sequential approximation
technique was employed by \cite{bdk05} in solving a Poisson
disorder problem in which the disorder distribution was
exponential and post disorder arrival rate was a random variable.
The authors formulated the problem under an auxiliary probability measure
as an optimal stopping time of an $\mathbb{R}_+$-valued \emph{odds-ratio} 
process. If we used a formulation similar to theirs, we would
obtain an optimal stopping problem with an unbounded continuation
region. Therefore, that formulation is not suitable for numerical
implementation. Also, the optimal stopping problem we consider
involves a terminal penalty term and a running cost with no
discount factor. In this section, we also show that an optimal
stopping time exists, and we describe two different types of
$\eps$-optimal stopping times. In Section~~\ref{sec:solution}, we
describe a numerical algorithm that can approximate the optimal
Bayes risk to a given level of accuracy. Finally,
Section~\ref{sec:examples} provides several examples illustrating
our solution. Appendix is home for the longer proofs.

\section{Problem Statement}\label{sec:prob-stat}
Let $(\Omega, \mathcal{F}, \mathbb{P})$ be a probability space hosting two
independent Poisson processes $(X_t^{(0)})_{t \geq 0}$, and
$(X_t^{(1)})_{t \geq 0}$ with intensities $\lambda_0$ and $\lambda_1$ respectively, and an
independent continuous-time Markov chain
$M=(M_t)_{t \geq 0}$ 
with state space
\begin{equation}
E \triangleq \{1,2,\cdots,n,\Delta\}.
\end{equation}
Here, $\Delta$ is an absorbing state and all the other states are
transient. The infinitesimal generator of $M$, which we denote by
$\mathcal{A}=(q_{i j})_{i,j \in E}$, is of the form
\begin{equation}\label{eq:matrix}
\mathcal{A}=
\begin{pmatrix}
&  & &
\\ &  R & & r
\\ & & & &
\\ 0 & \cdots &0 & 0
\end{pmatrix}
\end{equation}
where the $n \times 1$ vector $r$ is non-negative, and the $n \times n$ matrix $R$ is nonsingular.
The matrix $R$ has negative diagonal and nonnegative off-diagonal entries. Moreover $R$ and $r$
satisfy $R\cdot \vec{1} + r = \vec{0}$. 

For a point $\vec{\pi}=[\pi_1, \pi_2, \cdots, \pi_n, \pi] $ in
\begin{equation}
D\triangleq \{\vec{\pi}\in [0,1]^{n+1}: \sum_{i=1}^n \pi_{i} + \pi
=1\},
\end{equation}
let $\P$ denote the probability measure $\mathbb{P}$ such that the process $M$ has initial distribution $\vp$. That is,
\begin{equation}
\label{def:P-pi}
\P \{ A \} = \pi_1 \, \mathbb{P} \{A | M_0=1\} + \ldots
+ \pi_n \, \mathbb{P} \{A | M_0=n\}  + \pi \mathbb{P} \{A | M_0=\Delta\},
\end{equation}
for all $A \in \mathcal{F}$.
%
The absorption time of $M$ is defined as $\Theta \triangleq
\inf\{t>0: M_t =\Delta\}$, and its distribution is denoted by
\begin{equation}
\label{eq:distribution-of-theta-under-P} F_{\vp}(t) \triangleq
\P\{\Theta \leq t \}=1-[\pi_1, \pi_2, \cdots, \pi_n] \cdot
\exp(tR) \cdot \vec{1}, \quad 0 \leq t < \infty.
\end{equation}
Here, $\Theta$ is said to have a phase-type distribution, see e.g.
\cite{neuts}. 

The processes $X^0$, $X^1$ and $M$ are
unobservable. Rather we observe 
\begin{equation}
X_t=\int_0^{t}1_{\{s < \theta\}}dX_t^{(0)}+\int_0^{t}1_{\{s \geq
\theta\}}dX_t^{(1)},
\end{equation}
whose natural filtration will be denoted by $\Fb = \{\F_t\}_{t\ge
0}$. Let us define $\mathbb{G}=\{\mathcal{G}_t\}_{t \geq 0}$ as
an initial enlargement of $\mathbb{F}$ by setting
$\mathcal{G}_t \triangleq \mathcal{F}_t \vee \sigma\{M_t: t\geq 0
\}$. That is; $\mathcal{G}_t$ is the information available to a
\emph{genie} at time $t$ who is given the paths of the process
$M_t, t \geq 0$. If the paths of $M_t$, $t \geq 0$, are available
at time 0, then the observations come from a process $X$ that is a
Poisson process with rate $\lambda_0$ on the time interval
$[0,\Theta)$ and with rate $\lambda_1$ on $[\Theta,\infty)$ for
known positive constants $\lambda_0$ and $\lambda_1$.
Specifically, the observation process $X$ is a counting process such that
\text{$X_t - \int^t_{0} \left[\lambda_0 1_{\{s<\Theta\}}+
\lambda_1 1_{\{s\ge
    \Theta\}} \right]ds,\; t\ge 0$ is a
    $(\P,\mathbb{G})$-martingale}.
 The crucial feature here
is that $\Theta$ is neither known nor observable; only the process
$X$ is observable.  The problem is then to find a quickest
detection rule for the disorder time $\Theta$, which is
\emph{adapted to} the history $\Fb$ generated by the observed
process $X$ only. 
A detection rule is a stopping
time $\tau$ of the filtration $\mathbb{F}$, and we will denote the set
of these stopping times by $\mathcal{S}$. Our objective is to
find an element of $\mathcal{S}$ minimizing 
the Bayes risk
\begin{equation}\label{eq:bayes-risk}
 R_\tau(\vec{\pi}) \triangleq
\P\{\tau<\Theta\} +
  c\,\E(\tau-\Theta)^+, \quad
\end{equation}
for some positive constant $c$. Here $a^{+}=\max(a,0)$ for any $a
\in \R$. The first term in (\ref{eq:bayes-risk}) penalizes the
frequency of false alarms and the second term penalizes the
detection delay.
\begin{remark}
In order to minimize $R_\tau(\vec{\pi})$ in
$\mathcal{S}$, it is enough to consider stopping times with
bounded expectation. Indeed, if $\E\{\tau\}>1/c+E\{\Theta\}$, then
$R_{\tau}(\vec{\pi})\geq c (E\{\tau\}-E\{\Theta\} )> 1$, which is
greater than the cost incurred 
upon stopping immediately.
In the remainder we will use $\mathcal{S}_f$ to denote the class of $\mathbb{F}$-stopping
times whose expectation are strictly less than or equal to
$1/c+E\{\Theta\}$.
\end{remark}

Our objective is then to compute
\begin{equation}\label{eq:value-function}
V(\vec{\pi}) \triangleq \inf_{\tau \in \mathcal{S}_f
}R_{\tau}(\vec{\pi})=R_{\tau^*}(\vec{\pi}), \quad \text{for all}
\quad \vec{\pi} \in D,
\end{equation}
and to identify a rule $\tau^{*}$ (if there exists one) for which this infimum is attained. Note also that we have $0 \leq V(\vec{\pi}) \leq 1$ for all $\vec{\pi} \in D$.

\begin{remark}\label{rem:val-func}
Let us introduce the posterior probability distribution
$\vec{\Pi}_t \triangleq [\Pi^{(1)}_t, \cdots \Pi^{(n)}_t,\Pi_t]$, $t \geq 0$,
where
\begin{equation}\label{eq:posterior-prob}
\Pi_t \triangleq \P\{\theta \leq t| \F_t\}=\P\left\{M_t=\Delta|\F_t\right\}, \quad \text{and} \quad
\Pi^{(i)}_t \triangleq \P\left\{M_t=i|\F_t\right\},\quad t \geq 0,
\end{equation}
for $i \in \{1,\cdots,n\}$. Using the identities $\P\{\tau<\theta\}=\E\{1-\Pi_{\tau}\}$ and
\begin{equation}\label{eq:peanlty-in-t-Pi}
 \E\{(\tau-\Theta)^+\}=\E\left\{\int_0^{\tau}1_{\{\Theta\leq t
\}}dt\right\}=\E\left\{\int_0^{\infty}1_{\{\Theta\leq t \}}
1_{\{\tau>t\}} dt\right\}=\E\left\{\int_0^{\tau}\Pi_t dt\right\},
\end{equation}
we can represent the value function in (\ref{eq:value-function})
in terms of posterior probability distribution as
\begin{equation}\label{eq:value-func}
V(\vec{\pi})=\inf_{\tau \in
\mathcal{S}_f}\E\left\{\int_0^{\tau}k(\vec{\Pi}_t)dt+h(\vec{\Pi}_{\tau})\right\},
\end{equation}
in which
\begin{equation}
k(\vec{\pi}) \triangleq c \pi, \quad \text{and} \quad h(\vec{\pi})
\triangleq  1-\pi,
\end{equation}
 for all
$\vec{\pi}=[\pi_1, \pi_2, \cdots, \pi_n, \pi] \in D$.
\end{remark}

\begin{remark}
It follows from (\ref{eq:value-func}) that
\begin{equation}
V(\vec{\pi})\leq h(\vec{\pi})=1-\pi,
\end{equation}
for all $\vec{\pi}=[\pi_1, \pi_2, \cdots, \pi_n, \pi] \in D$.
\end{remark}

\begin{lemma}
Let us define the hazard rate of the distribution $F_{\vp}$ of $\Theta$ as
\begin{equation}
\eta(t)=\frac{F_{\vp}'(t)}{1-F_{\vp}(t)}, \qquad \text{ for $t > 0$}.
\end{equation}
Then the a-posteriori probability process $(\Pi_t)_{t \geq 0}$
satisfies
\begin{equation}\label{eq:pi-h}
d\Pi_t=[\eta(t)-(\lambda_1-\lambda_0)\Pi_t](1-\Pi_t)dt+\frac{(\lambda_1-\lambda_0)\Pi_{t-}
(1-\Pi_{t-}) }{\lambda_0 (1-\Pi_{t-})+\lambda_1 \Pi_{t-}}dX_t,
\end{equation}
with $\Pi_0=\pi$. 
\end{lemma}
\begin{proof}
We will first introduce a reference probability measure $\P_0$
under which the processes $M$ and $X$ are independent. Moreover, the
probability law of $M$ under $\P_0$ will remain unchanged.

Let us introduce
\begin{align}
  Z_t \triangleq \exp \left\{\int^t_0 \log
    \left(\frac{H(s)}{\lambda_0}\right)\, dX_s - \int^t_0
    [H(s)-\lambda_0] ds \right\}, \quad t\ge 0,
\end{align}
in which $H(s) \triangleq \lambda_0 1_{\{s<\Theta\}}+ \lambda_1
1_{\{s\ge
    \Theta\}}$.
 Using the process $Z$ we can define a new probability measure $\P_0$ on
$(\Omega,  \mathbb{G})$ locally in terms of the Radon-Nikodym derivatives 
\begin{align}
\label{eq:Radon-Nikodym-derivative}
  \left.\frac{d\P}{d\P_0} \right|_{\G_t} = \frac{1}{Z_t} =
  1_{\{\theta>t\}} + 1_{\{\theta \le t\}} \frac{L_{\theta}}{L_t}
\end{align}
for every $0\le t <\infty$, where
\begin{align}
\label{eq:likelihood-ratio-when-both-rates-are-known}
  L_t \triangleq \left(\frac{\lambda_1}{\lambda_0}
  \right)^{X_t}e^{-(\lambda_1-\lambda_0)t}.
\end{align}
Under the measure $\P_0$, the process $Z$ is a martingale, $X$ is a Poisson process with
intensity $\lambda_0$ and is independent of $M$ (see e.g.
Section 2 in \cite{bdk05}, or Appendix A1 in \cite{ds}). Moreover,
$\P$ and $\P_0$ coincide on $\mathcal{G}_0=\sigma\{ M_s; s\ge 0 \} $,
therefore $\P_0\{\Theta \leq t\}=F_{\vp}(t)$.

Using the Bayes rule (see e.g. \cite{MR2001k:60001a}) (this is
also known as the Kallianpur-Striebel formula) we obtain
\begin{equation}\label{eq:piandbarpi}
\Pi_t=\P\{\Theta \leq t|\F_t\}=\frac{\E_0\{Z_t 1_{\{\Theta \leq t
\}}| \F_t \}}{\E_0\{Z_t|\F_t\}}, \quad
1-\Pi_t=\frac{\E_0\{1_{\{\Theta>t\}}|\F_t\}}{\E_0\{Z_t|\F_t\}}=\frac{(1-F_{\vp}(t))}{\E_0\{Z_t|\F_t\}}.
\end{equation}
Here, to derive the second equality in the second equation we used
independence of $\theta$ and $X$ under $\P_0$.

Let us define \emph{the odds ratio} process
\begin{equation}\label{eq:odds-ratio}
\Phi_t \triangleq \frac{\Pi_t}{1-\Pi_t}, \quad 
0 \leq t <\infty.
\end{equation}
Using (\ref{eq:piandbarpi}) we can obtain a new representation for
the odds ratio process
\begin{equation}\label{eq:phi-t-L}
\Phi_t=\frac{\E_0\{Z_t 1_{\{\Theta \leq t \}}| \F_t
\}}{(1-F_{\vp}(t))} =\frac{1}{1-F_{\vp}(t)}\left[\pi
L_t + \int_0^t \frac{L_t}{L_s} F'_{\vp}(s) ds\right],
\end{equation}
Here, again we used the independence of $\Theta$ and $\mathbb{F}$.
The process $L=\{L_t, t\ge 0\}$ in
(\ref{eq:Radon-Nikodym-derivative}) is a $(\P_0, \Fb)$-martingale
and is the unique locally bounded solution of the equation
\begin{align*}
dL_t = [(\lambda_1/\lambda_0)-1]L_{{t-}} (dX_t-\lambda_0
dt),\qquad L_0=1;
\end{align*}
see, e.g.,  \cite{RY99} or \cite{MR2003j:60001}. Applying chain
rule to (\ref{eq:phi-t-L}), we get
\begin{equation}\label{eq:dyn-phi-with-h}
d \Phi_t= \eta(t)(1+\Phi_t)dt+
\Phi_{t-}\left(\frac{\lambda_1}{\lambda_0}-1\right)d(X_t-\lambda_0
t), \quad \Phi_0=\frac{\pi}{1-\pi}.
\end{equation}
By an another application of chain rule to (\ref{eq:odds-ratio})
together with (\ref{eq:dyn-phi-with-h}) we obtain (\ref{eq:pi-h}).
\end{proof}

\begin{proposition} \label{prop:dynamics-of-pi}
The dynamics of the posterior probability distribution
$\vec{\Pi}_t=[\Pi^{(1)}_t, \cdots \Pi^{(n)}_t,\Pi_t]$, $t \geq 0$,
which is defined in (\ref{eq:posterior-prob}), is given by
\begin{align}
d\Pi_t&=\left(\sum_{j=1}^{n}q_{j
\Delta}\Pi^{(j)}_t-(\lambda_1-\lambda_0)\Pi_t(1-\Pi_t)\right)dt+\frac{(\lambda_1-\lambda_0)\Pi_{t-}(1-\Pi_{t-})}
{\lambda_0(1-\Pi_{t-})+\lambda_1 \Pi_{t-}}dX_t, 
\label{eq:dyn-pi-t}
\\ d\Pi_t^{(i)}&=\left(\sum_{j=1}^{n}q_{j i}\Pi^{(j)}+(\lambda_1-\lambda_0)\Pi_t
\Pi_t^{(i)}\right)dt-
\frac{(\lambda_1-\lambda_0)\Pi_{t-}\Pi^{(i)}_{t-}}
{\lambda_0(1-\Pi_{t-})+\lambda_1 \Pi_{t-}}dX_t,
\label{eq:dyn-pi-i-t}
\end{align}
for $i \in \{1,\cdots,n\}$, and with $\vec{\Pi}_0=[\pi_{1}, \cdots \pi_{n},\pi]$.
\end{proposition}

\begin{proof}
First, observe that the hazard rate function of the distribution
$F_{\vp}$, can be written as
\begin{equation}\label{eq:hr}
\eta(t)=\frac{\sum_{i=1}^{n} \P\{M_t=i\} q_{i
\Delta}}{1-F_{\vp}(t)}.
\end{equation}
On the other hand,
\begin{align}\label{eq:pi-t}
\begin{aligned}
\Pi_t^{(i)}=\E\{1_{\{M_t=i\}}|\F_t\} &=\frac{\E_0\{Z_t
1_{\{M_t=i\}}|\F_t\}}{\E_0\{Z_t|\F_t\}} \\ &=\frac{\E_0\{
1_{\{M_t=i\}}|\F_t\}}{\E_0\{Z_t|\F_t\}}=\frac{\E_0\{
1_{\{M_t=i\}}\}}{\E_0\{Z_t|\F_t\}}=\frac{\E\{
1_{\{M_t=i\}}\}}{\E_0\{Z_t|\F_t\}},
\end{aligned}
\end{align}
in which $\E_0$ denotes the expectation under the measure $\P_0$
which we introduced in (\ref{eq:Radon-Nikodym-derivative}). The
second equality in this equation follows from Bayes' formula, the
third equality follows from the definition of $Z$ in
(\ref{eq:Radon-Nikodym-derivative}), the fourth equality follows
from the independence of $M$ and $X$ under the measure $\P_0$,
and, finally, the fourth equality follows from the fact that under
the measures $\P$ and $\P_0$ the law of $M$ is the same.

From (\ref{eq:piandbarpi}) and (\ref{eq:pi-t}) it is immediate
that
\begin{equation}\label{eq:pi-i-over-1-pi}
\frac{\Pi^{(i)}_t}{1-\Pi_t}=\frac{\P\{M_t=i\}}{1-F_{\vp}(t)},
\quad i\in\{1, \cdots,n\}.
\end{equation}
Then, from (\ref{eq:pi-t}) and (\ref{eq:pi-i-over-1-pi}) it follows
that
\begin{equation}\label{eq:alt-for-h}
\eta(t)=\frac{\sum_{i=1}^{n}q_{i \Delta} \Pi_t^{(i)}}{1-\Pi_t}.
\end{equation}
This equation together with (\ref{eq:pi-h}) yields
(\ref{eq:dyn-pi-t}).

We will now derive the dynamics of $(\Pi_t^{(i)})_{t \geq 0}$, $i
\in \{1,\cdots n \}$. Let $p_{ij}(t) \triangleq \P\{M_t=j|M_0=i\}$ denote the transition probabilities of the process $M$. Recall that $t \rightarrow p_{ij}(t)$, $t
\geq 0$, satisfies the forward Kolmogorov equation, i.e.,
\begin{equation}\label{eq:kolmogorov}
\frac{dp_{ij}(t)}{dt}=\sum_{k=1}^{n}q_{k j}p_{ik}(t).
\end{equation}
and that
\begin{equation}\label{eq:markov-chain-prob-m-i}
\P\{M_t=i\}=\sum_{j=1}^{n}\pi_j p_{ji}(t).
\end{equation}
Now, applying chain rule to (\ref{eq:pi-i-over-1-pi}) we obtain
\begin{equation}
\begin{split}
d\Pi_t^{(i)}&=-\frac{\Pi_t^{(i)}}{1-\Pi_t} d\Pi_t+
(1-\Pi_t)\frac{\sum_{j=1}^{n}\pi_j \sum_{k=1}^{n} q_{ki}
p_{jk}(t)+\sum_{j=1}^{n} \pi_j p_{ji}(t) \eta(t)}{1-F_{\vp}(t)} dt
\\&= -\frac{\Pi_t^{(i)}}{1-\Pi_t} d\Pi_t+(1-\Pi_t) \left(
\frac{\sum_{k=1}^{n} \P(M_t=k) q_{ki}}{1-F_{\vp}(t)}+
\eta(t)\frac{\P\{M_t=i\}}{1-F_{\vp}(t)} \right) dt
\\&=-\frac{\Pi_t^{(i)}}{1-\Pi_t} d\Pi_t+ \left( \sum_{k=1}^{n}\Pi_{t}^{(k)}q_{ki}+\eta(t)\Pi^{(i)}_t \right) dt.
\end{split}
\end{equation}
The first line follows from (\ref{eq:kolmogorov}), and the second
follows from (\ref{eq:markov-chain-prob-m-i}). The last line is a
result of the identity in (\ref{eq:pi-i-over-1-pi}). This
equation, together with (\ref{eq:dyn-pi-t}) and (\ref{eq:alt-for-h}) gives
(\ref{eq:dyn-pi-i-t}).
\end{proof}

\begin{remark} \label{rem:deterministic-paths}
Let $\vx(t,\vp) \triangleq (x_1(t,\vp), \cdots , x_{n}(t,\vp),
x_{\Delta}(t,\vp))$ be the solution of the system of ordinary
differential equations
\begin{equation}\label{eq:dyn-x}
\begin{split}
dx_{\Delta}(t, \vp)&=\left(\sum_{j=1}^{n}q_{j \Delta} x_{j}(t,
\vp)- (\lambda_1-\lambda_0)x_{\Delta}(t, \vp)(1-x_{\Delta}(t,
\vp)) \right)dt, \quad \text{with }x_{\Delta}(0 ,\vp)=\pi,
\\
 dx_{i}(t,\vp)&=\left(\sum_{j=1}^{n} q_{ji} x_j(t, \vp)+
(\lambda_1-\lambda_0) x_{\Delta}(t,\vp) x_{i}(t)\right)dt, \quad
\text{with }x_{i}(0 , \vp)=\pi_i, 
\end{split}
\end{equation}
for $i \in \{1, \cdots,n\}$. Due to Kolmogorov's forward equations,
the solution of this system of equations can be written as
\begin{align}
\label{eq:solns-of-system}
\begin{aligned}
x_{\Delta}(t,\vp)&=  \frac{\pi e^{-(\lambda_1-\lambda_0)t}+ 
\int_{0}^{t}e^{-(\lambda_1-\lambda_0)(t-s)}F'_{\vp}(s)ds}{1-F_{\vp}(t)+
\pi e^{-(\lambda_1-\lambda_0)t}+ 
\int_{0}^{t}e^{-(\lambda_1-\lambda_0)(t-s)}F'{\vp}_(s)ds},
\\
x_{i}(t,\vp)&= \frac{\sum_{j=1}^{n} \pi_j
p_{ji}(t)}{1-F_{\vp}(t)+ \pi e^{-(\lambda_1-\lambda_0)t}+ 
\int_{0}^{t}e^{-(\lambda_1-\lambda_0)(t-s)}F'_{\vp}(s)ds}, \qquad \text{for }i \in
\{1, \cdots n\},
\end{aligned}
\end{align}
in terms of the transition probabilities
$p_{ij}(t) \triangleq \P\{M_t=j|M_0=i\}$, for $i,j \in E$.
Moreover, the expressions in (\ref{eq:solns-of-system}) are equivalent to
\begin{equation}\label{eq:alternative-esp-for-x}
\begin{split}
x_{\Delta}(t,
\vp)=\frac{\E\left\{1_{\{t \geq
\theta\}}e^{-(\lambda_1-\lambda_0)(t-\theta)}\right\}}{\E\left\{e^{-(\lambda_1-\lambda_0)(t-\theta)^+}\right\}}, \quad \text{and} \quad
x_{i}(t,\vp)=\frac{\P\left\{M_t=i\right\}}{\E\left\{e^{-(\lambda_1-\lambda_0)(t-\theta)^{+}}\right\}},
\end{split}
\end{equation}
for $i \in \{1,\cdots,n\}$.
\end{remark}

Using Markov property of $M$ and (\ref{eq:alternative-esp-for-x}), we have
\begin{align}
\label{eq:probability-for-ts}
\begin{aligned}
\P \{ M_{t+s} = i\} &= \sum_{ j=1}^n \P\{ M_t =j\} \cdot  \P \{ M_{t+s} = i | M_t =j \} \\
&=   \E\left\{e^{-(\lambda_1-\lambda_0)(t-\theta)^+}\right\} \sum_{ j=1}^n x_j (t, \vp ) \cdot  \P \{ M_{t+s} = i | M_t =j \}
\\&
=  \E\left\{e^{-(\lambda_1-\lambda_0)(t-\theta)^+}\right\} \cdot \mathbb{P}^{\vx(t,\vp)} \{ M_{s} = i\}
\end{aligned}
\end{align}
for $i \le n$, and 
\begin{align}
\label{eq:expectation-for-ts}
\begin{aligned}
 \E &\left\{e^{-(\lambda_1-\lambda_0)(t+s-\theta)^+}\right\} = 
  \E\left\{  \E \left\{ e^{-(\lambda_1-\lambda_0)(t+s-\theta)^+} \Big|  M_s ; s\le t \right\} \right\} \\
  &= \E \left\{\sum_{j=1}^n 1_{ \{ M_t =j \} } \cdot \E \left\{ e^{-(\lambda_1-\lambda_0)(s-\theta)^+} \big|  M_0 =j  \right\}
  + 1_{ \{ M_t =\Delta \} } \cdot e^{-(\lambda_1-\lambda_0)(t+s-\theta)}  \right\} \\
 &= \sum_{j=1}^n \P \{ M_t =j \} \E \left\{ e^{-(\lambda_1-\lambda_0)(s-\theta)^+} \big|  M_0 =j  \right\}
 +  e^{-(\lambda_1-\lambda_0)s} \cdot \E \left\{ 1_{ \{ t \ge \Theta \} } \cdot e^{-(\lambda_1-\lambda_0)(t-\theta)}  \right\} \\
 &=  \E\left\{e^{-(\lambda_1-\lambda_0)(t-\theta)^+}\right\}  \left[ 
 \sum_{j=1}^n x_{j} (t, \vp)\cdot  \E \left\{ e^{-(\lambda_1-\lambda_0)(s-\theta)^+} \big|  M_0 =j  \right\}
 + x_{\Delta}(t,\vp) \cdot e^{-(\lambda_1-\lambda_0)s} 
 \right] \\
 &=  \E\left\{e^{-(\lambda_1-\lambda_0)(t-\theta)^+}\right\}  \cdot  
 \mathbb{E}^{\vx(t,\vp)} \left\{e^{-(\lambda_1-\lambda_0)(s-\theta)^+}\right\} 
\end{aligned}
\end{align}

Using (\ref{eq:probability-for-ts}) and (\ref{eq:expectation-for-ts}), it is now easy to see that $t \mapsto x(t,\vp)$ has the semi-group property $x(t+s, \vp ) = x( t, x(s,\vp) ) $. Then, the dynamics in (\ref{eq:dyn-pi-t}), (\ref{eq:dyn-pi-i-t}) and
Remark~\ref{rem:deterministic-paths}
imply that $\vP$ is a piecewise deterministic 
process whose natural filtration coincides with $\mathbb{F}$. Between two jumps of $X$,
the process $\vP$ follows the curves $t \mapsto \vx(t,\vp) $, and at arrival times of $X$, it jumps from one curve to another. More precisely, the paths of $\vP$ have the characterization
     \begin{equation}\label{eq:rel-pi-x}
     \begin{split}
\vP_t=\vx\left(t-\sigma_m,\vP_{\sigma_m}\right), \quad
\quad \sigma_m \leq t< \sigma_{m+1}, \;\; m\in \mathbb{N}_0, \qquad\qquad\qquad\qquad \\
\vP_{\sigma_m}=\left(\frac{\lambda_0
\Pi^{(1)}_{\sigma_m-}}{\lambda_0(1-\Pi_{\sigma_m-})+\lambda_1
\Pi_{\sigma_m-} }, \cdots, \frac{\lambda_0
\Pi^{(b)}_{\sigma_m-}}{\lambda_0(1-\Pi_{\sigma_m-})+\lambda_1
\Pi_{\sigma_m-} }, \frac{\lambda_1
\Pi_{\sigma_m-}}{\lambda_0(1-\Pi_{\sigma_m-})+\lambda_1
\Pi_{\sigma_m-} }\right),
\end{split}
\end{equation}
%
in which
\begin{equation}\label{eq:sigma}
\sigma_0\equiv 0, \quad \text{and} \quad \sigma_m \triangleq
\inf\{t>\sigma_{m-1}|X_t-X_{t-}>0\}, \qquad m \in \N.
\end{equation}
Moreover, for a bounded function $g(\cdot)$, we have 
\begin{align}
\label{eq:Markov-justification}
\begin{aligned}
\E &\left\{ g( X_{t+s} - X_t ) \big| \F_t \right\}\\ 
&= \sum_{j=1}^n \P \{ M_t =j \big| \F_t  \} 
\cdot \E \left\{ g( X_{t+s} - X_t ) \big| \F_t , M_t =j\right\}  
\\ &\qquad \qquad \qquad \qquad \qquad \qquad \qquad 
+\P \{ M_t = \Delta \big| \F_t  \} \E \left\{ g( X_{t+s} - X_t ) \big| \F_t , M_t =\Delta \right\} \\
&= \sum_{j=1}^n \Pi^{(i)}_t
\cdot \E \left\{ g( X_{s}  ) \big|  M_0 =j\right\}   
+ \Pi_t \cdot \E \left\{ g( X_{s})  \big| M_0 =\Delta \right\} = \mathbb{E}^{\vP_t}  \left\{ g( X_{s} )\right\} .
\end{aligned}
\end{align}
Then, the characterization in (\ref{eq:rel-pi-x}) and (\ref{eq:sigma}) implies 
that $\vP$ is a $(\P, \Fb)$-Markov process due to (\ref{eq:Markov-justification}). 

%


\section{Sequential Approximation}\label{sec:sequential}
Let us define the sequence of functions
\begin{equation}\label{eq:auxiliary}
V_{m}(\vec{\pi})=\inf_{\tau \in
\mathcal{S}_f}\E\left\{\int_{0}^{\tau \wedge \sigma_m}
k(\vec{\Pi}_t)dt+ h(\vec{\Pi}_{\tau \wedge \sigma_m})\right\},
\end{equation}
in which $\sigma_m$, $m \in \mathbb{N}_0$, is defined in
(\ref{eq:sigma}). The functions $V_m(\cdot)$, $m \in
\mathbb{N}_0$, are non-negative and bounded above by
$h(\cdot)$. By definition, the sequence $\{V_m\}_{m \geq 1}$ is
decreasing and $V_m \geq V$ for all $m$. Therefore the point-wise
limit $\lim_{m \rightarrow \infty} V_m$ exists and is greater than
or equal to $V$. In fact a stronger convergence result holds as
the next lemma shows.

\begin{proposition}\label{prop:sequential}
As $m \rightarrow \infty$, the sequence $\{V_m(\cdot)\}_{m \geq
1}$ converges to $V(\cdot)$ uniformly on $D$. In fact, for every
$m \in \mathbb{N}$
\begin{equation}\label{eq:uniform-conv}
V_m(\vec{\pi})-\sqrt{\left(\frac{1}{c}+\E\{\theta\}\right)\frac{\max\{\lambda_0,
\lambda_1\}}{m-1}}\leq V(\vec{\pi}) \leq V_m(\vec{\pi}), \quad
\text{for all} \quad \vec{\pi} \in [0,1]^{n+1}.
\end{equation}
\end{proposition}

\begin{proof} The second inequality in (\ref{eq:uniform-conv}) follows immediately, since by definition $V_m(\cdot) \geq V(\cdot)$.
Let us prove the first inequality. For any $\tau \in
\mathcal{S}_f$, the expectation
$\E\left\{\int_0^{\tau}k(\vec{\Pi}_t)dt+
h(\vec{\Pi}_{\tau})\right\}$ can be written as
\begin{equation}\label{eq:ineq}
\begin{split}
& \E\left\{\int_0^{\tau \wedge \sigma_m }k(\vec{\Pi}_t)dt+
h(\vec{\Pi}_{\tau \wedge \sigma_m})\right\}
+\E\left\{1_{\{\tau>\sigma_m\}}
\left[\int_{\sigma_m}^{\tau}k(\vec{\Pi}_t)dt+h(\vec{\Pi}_{\tau})-h(\vec{\Pi}_{\sigma_m})\right]\right\}
\\&  \geq \E\left\{\int_0^{\tau \wedge \sigma_m
}k(\vec{\Pi}_t)dt+h(\vec{\Pi}_{\tau \wedge
\sigma_m})\right\}-\E\left\{1_{\{\tau>\sigma_m\}}\right\},
\end{split}
\end{equation}
since $0 \le h(\cdot) \le 1$. Note that
\begin{equation}\label{eq:cauchy}
\E\left\{1_{\{\tau>\sigma_m\}}\right\}\leq
\E\left\{1_{\{\tau>\sigma_m\}}\left(\frac{\tau}{\sigma_m}\right)^{1/2}\right\}
\leq \E\left\{\left(\frac{\tau}{\sigma_m}\right)^{1/2}\right\}
\leq \sqrt{\E\{\tau\} \E\left\{\frac{1}{\sigma_m}\right\}},
\end{equation}
which follows as a result of Cauchy-Schwartz inequality, and that
\begin{equation}\label{eq:inverse-moment}
 \E\left\{\frac{1}{\sigma_m}\right\} \leq \frac{\max\{\lambda_0,\lambda_1\}}{m-1}.
\end{equation}
Since $\E\{\tau\} \leq 1/c+\E\{\theta\}$ for any $\tau \in
\mathcal{S}_f$, using (\ref{eq:ineq}), (\ref{eq:cauchy}) and
(\ref{eq:inverse-moment}) we obtain
\begin{multline}
\E\left\{\int_0^{\tau}k(\vec{\Pi}_t)dt+
h(\vec{\Pi}_{\tau})\right\} \leq \\  \E\left\{\int_0^{\tau \wedge
\sigma_m }k(\vec{\Pi}_t)dt+ h(\vec{\Pi}_{\tau \wedge
\sigma_m})\right\}-\sqrt{\left(\frac{1}{c}+\E\{\theta\}\right)\frac{\max\{\lambda_0,
\lambda_1\}}{m-1}}.
\end{multline}
Now taking the infimum of both sides over the stopping rules in
$\mathcal{S}_f$, we obtain the first inequality in
(\ref{eq:uniform-conv}).
\end{proof}

To calculate the functions $V_m(\cdot)$ iteratively, we introduce
the following operators acting on bounded functions $w:D
\rightarrow \R$
\begin{equation}\label{eq:defn-J}
\begin{split}
J w(t, \vec{\pi})&=\E\left\{\int_0^{t \wedge
\sigma_1}k(\vec{\Pi}_s)ds+1_{\{t<\sigma_1\}}h(\vec{\Pi}_t) +1_{\{t
\geq \sigma_1\}}w(\vec{\Pi}_{\sigma_1})\right\} \quad t \in
[0,\infty],
\\  J_t w(\vec{\pi})&=\inf_{u \in [t,\infty]}J w(u,\vec{\pi}), \quad t \in [0,\infty].
\end{split}
\end{equation}
The action of the operator $J$ on the function $w$ can be written
as
\begin{equation}\label{eq:exp-for-J-w}
Jw(t,\vec{\pi})=\int_0^{t}\P\{s \leq
\sigma_1\}k(\vec{x}(s,\vec{\pi}))ds+\int_0^{t}\P\{\sigma_1 \in
ds\} Sw(\vec{x}(s,\vec{\pi}))+h(\vec{x}(t,\pi))\P\{t<\sigma_1\},
\end{equation}
in which
\begin{equation}
Sw(\vec{\pi})\triangleq w\left(\frac{\lambda_0
\pi_1}{\lambda_0(1-\pi)+\lambda_1 \pi}, \cdots, \frac{\lambda_0
\pi_n}{\lambda_0(1-\pi)+\lambda_1 \pi}, \frac{\lambda_1
\pi}{\lambda_0(1-\pi)+\lambda_1 \pi}\right).
\end{equation}
Let us now compute the distribution and the density of $\sigma_1$
under $\P$, respectively, since it appears in the expression for
$Jw$. We have
\begin{equation}\label{eq:distribution-of-first-jump}
\begin{split}
\P\{\sigma_1>t\}&=\int_0^{\infty}\P\{\sigma_1>t|\theta \in ds\}
\P\{\theta \in ds\}
\\&=\int_0^{t}\P\{\sigma_1>t|\theta=s\}\P\{\theta \in
ds\}+\int_t^{\infty}\P\{\sigma_1>t|\theta=s\}\P\{\theta \in ds\}
\\&=\int_0^{t}e^{-\lambda_0 s}
e^{-\lambda_1(t-s)}\P\{\theta \in ds
\}+\int_t^{\infty}e^{-\lambda_0 t}\P\{\theta \in ds \}
\\&=e^{-\lambda_0 t}
\E\left\{e^{-(\lambda_1-\lambda_0)(t-\theta)^{+}}\right\},
\end{split}
\end{equation}
from which it follows that
\begin{equation}\label{eq:sigma-density}
\begin{split}
\P\{\sigma_1 \in dt\}
=e^{-\lambda_0 t}\left[\lambda_0 \cdot \E\{1_{\{t<\theta\}}\}+\lambda_1\cdot
\E\left\{1_{\{t \geq
\theta\}}\, e^{-(\lambda_1-\lambda_0)(t-\theta)}\right\}\right]dt.
\end{split}
\end{equation}

\begin{remark}
\label{rem:inf-J-is-attained}
For a bounded function $w(\cdot)$,
using equations in 
(\ref{eq:alternative-esp-for-x}), (\ref{eq:distribution-of-first-jump}) and (\ref{eq:sigma-density}), it can be verified easily that the integrands in
(\ref{eq:exp-for-J-w}) are absolutely integrable. Hence
\begin{align*}
\lim_{t \to \infty } Jw(t, \vp) = Jw(\infty,\vp) < \infty,
\end{align*}
and the mapping $t \to Jw(t,\vp)$ is continuous on $[0,\infty]$. Therefore, the infimum
in (\ref{eq:defn-J}) is attained for all $t \in [0,\infty]$.
\end{remark}

\begin{remark}\label{rem:monotone}
\begin{itemize}
\item[(i)] $0 \leq J_0 w(\cdot) \leq h(\cdot)$ for all non-negative and bounded
function $w$.
\item[(ii)] For two bounded functions $w_1(\cdot) \leq w_2(\cdot)$, we
have $J_0 w_1(\cdot) \leq J_0 w_2(\cdot)$.
\end{itemize}
\end{remark}

\begin{lemma}\label{lem:concave}
If $w:D\rightarrow \R$ is positive and concave, then so are the
mappings
\begin{equation}
\vec{\pi} \rightarrow J w(t,\vec{\pi}) \quad \text{and} \quad
\vec{\pi} \rightarrow J_0 w(\vec{\pi}).
\end{equation}
\end{lemma}

\begin{lemma}\label{lem:cont}
If the function $w:D \rightarrow \R_{+}$ is bounded and
 continuous, then $(t,\vec{\pi}) \rightarrow J_0 w
 (t,\vec{\pi})$ and $\vec{\pi} \rightarrow J_0 w(\vec{\pi})$ are
 also continuous functions.
\end{lemma}

Using the operator $J_0$, let us define a sequence of functions
\begin{equation}\label{eq:seq-of-func}
v_0(\vec{\pi})\equiv h(\vec{\pi}) \quad \text{and} \quad
v_m(\vec{\pi}) \triangleq J_0 v_{m-1}(\vec{\pi}), \, m \geq 1,
\quad \text{for all $\vp \in D$}.
\end{equation}

\begin{cor}\label{cor:little-v-n}
 Each $v_m(\cdot)$ is positive, continuous,
concave on $D$. The sequence $\{v_m(\cdot)\}_{m \geq 1}$ is
decreasing, hence the pointwise limit $v(\vec{\pi})=\lim_{m
\rightarrow \infty}v_{m}(\vec{\pi})$, $\vec{\pi} \in D$, exists.
The function $v(\cdot)$ is again concave.
\end{cor}
\begin{proof}
The proof easily follows from Remark~\ref{rem:monotone},
Lemmata~\ref{lem:concave} and \ref{lem:cont}. To prove the
concavity of $v(\cdot)$ we also use the fact that the lower
envelope of concave functions is concave.
\end{proof}

The following lemma, which follows from \cite{bremaud} Theorem
T.33, characterizes the stopping times of piece-wise deterministic
Markov processes. Also see \cite{MR96b:90002}, Theorem A2.3.
\begin{lemma}\label{lem:bremaud}
For every $\tau \in \mathcal{S}$, and for every $m \in
\mathbb{N}$, there exists a $\F_{\sigma_m}$-measurable random
variable such that $\tau \wedge \sigma_{m+1}=(\sigma_m+R_m) \wedge
R_{m+1}$, $\P$-almost surely on $\{\tau \geq \sigma_m\}$.
\end{lemma}

\begin{proposition}\label{prop:V-n-epsilon}
 For every
$\eps \geq 0$, let us define
\begin{equation}\label{eq:defn-r-m}
r_{m}^{\eps}(\vec{\pi}) \triangleq \inf\{s \in (0,\infty]: J
v_{m}(s,\vec{\pi}) \leq J_0 v_{m} (\vec{\pi})+\eps\}, \quad
\vec{\pi} \in D,
\end{equation}
\begin{equation}\label{eq:defn-S-eps}
S_1^{\eps} \triangleq r_0^{\eps}(\vec{\Pi}_0) \wedge \sigma_1
\quad \text{and} \quad S_{m+1}^{\eps}(\vec{\pi}) \triangleq
\begin{cases}
r_{m}^{\eps/2}(\vec{\Pi}_0) & \text{if
$\sigma_1>r_m^{\eps/2}(\vec{\Pi}_0)$},
\\ \sigma_1+ S_m^{\eps/2}\circ \theta_{\sigma_1} & \text{if
$\sigma_1 \leq r_m^{\eps/2}(\vec{\Pi}_0)$},
\end{cases}
\end{equation}
where $\theta_s$ is the shift operator on $\Omega$, i.e.,
$X_{t}\circ \theta_s=X_{s+t}$. Then, for every $m \geq 1$
\begin{equation}\label{eq:eps-opt}
\E\left\{\int_0^{S_m^{\eps}}k(\vec{\Pi}_t)dt+h(\vec{\Pi}_{S_m^{\eps}})\right\}
\leq v_{m}(\vec{\pi})+\eps.
\end{equation}
Moreover, for all $m \in \mathbb{N}$, $v_{m}(\vp)=V_m(\vp)$ on $D$.
\end{proposition}

\begin{proposition}\label{prop:v-V}
We have $v(\vp)=V(\vp)$ 
for every $\vec{\pi} \in D$.
Moreover, $V$ is the largest solution of $U=J_0 U$ that is smaller
than or equal to $h$.
\end{proposition}
\begin{lemma}\label{lem:dyn-p-J-t-J}
For every bounded function $\vp \rightarrow w(\vp)$, $\vp \in D$,
we have
\begin{equation}\label{eq:J-t-J}
J_t w(\vp)=J w (t,\vp)+\P\{\sigma_1>t\} \cdot \left\{ J_0 w
(\vx(t,\vp))-h(\vx(t,\vp))\right\},
\end{equation}
for all $t \geq 0$.
\end{lemma}

\begin{cor}\label{cor:r-m}
Let
\begin{equation}
r_{m}(\vp)=\inf\{t \in (0,\infty]:J v_{m}(s,\vp)=J_0 v_{m}(\vp)\}.
\end{equation}
Then
\begin{equation}
r_{m}(\vp)=\inf\{t \in (0,\infty]:
v_{m+1}(\vx(t,\vp))=h(\vx(s,\vp))\}.
\end{equation}
Here, we use the convention that $\inf \emptyset=\infty$.
\end{cor}

\begin{remark}
Substituting $w=v_m$ in (\ref{eq:J-t-J}) yields the dynamic
programming equation for the sequence of function $\{v_m
(\cdot)\}_{m \in \mathbb{N}_0}$; for every $\vp \in D$ and $n\in
\mathbb{N}_0$
\begin{equation}
v_{m+1}(\vp)=Jv_m(t,\vp)+\P\left\{\sigma_1>t\right\}  \cdot [v_{m+1}(\vx(t,\vp))-h(\vx(t,\vp))],
\quad \; t \in [0, r_m(\vp)].
\end{equation}
Moreover, if we take $w=V$ in (\ref{eq:J-t-J}), then we obtain
\begin{equation}\label{eq:dyn-prog-0.5}
J_t
V(\vp)=JV(t,\vp)+\P\left\{\sigma_1>t\right\} \cdot \left[V(\vx(t,\vp))-h(\vx(t,\vp))\right],
\quad t\geq 0.
\end{equation}
Let us define
\begin{equation}
r(\vp)\triangleq \inf\{t \in (0,\infty]:JV(t,\vp)=J_0 V(\vp)\}.
\end{equation}
The same arguments as in the proof of Corollary~\ref{cor:r-m}
leads to
\begin{equation}
r(\vp)=\inf\{t \in (0,\infty]: V(\vx(t,\vp))=h(\vx(t,\vp))\}.
\end{equation}
This equation together with (\ref{eq:dyn-prog-0.5}) yields
\begin{equation}
V(\vp)=JV(t,\vp)+\P\left\{\sigma_1>t\right\} \cdot [V(\vx(t,\vp))-h(\vx(t,\vp))],
\quad t \in [0,r(\vp)].
\end{equation}
\end{remark}

\begin{remark}\label{rem:right-cont-of-paths}
From Propositions~\ref{prop:sequential}, \ref{prop:V-n-epsilon},
\ref{prop:v-V} and Corollary~\ref{cor:little-v-n} it follows that
a continuous sequence of functions uniformly converge to $v=V$.
Therefore $V$ is continuous on $D$. Since $t \rightarrow
\vx(t,\vp)$, $t \geq 0$,  is continuous for all $\vp \in D$, the
mapping $t \rightarrow V(\vx(t,\vp))$, $t \geq 0$ is also
continuous for all $\vp \in D$. Moreover, for every $\vp$, the
path $t \rightarrow \vP_t$, $t \geq 0$, follows the deterministic
curves $t \rightarrow \vx(t,\vp)$ between the jumps. Hence the
process $t \rightarrow V(\vP_t)$ is right-continuous with left
limits.
\end{remark}

Let us define the $\mathbb{F}$-stopping times
\begin{equation}\label{eq:defn-U-eps}
U_{\eps} \triangleq \inf\{t \geq 0: V(\vP_t)-h(\vP_t) \geq
-\eps\}, \quad \eps \geq 0.
\end{equation}
Remark~\ref{rem:right-cont-of-paths} implies
\begin{equation}
V(\vP_{U_{\eps}})-h(\vP(U_{\eps})) \geq -\eps \quad \text{on the
event} \quad \{U_{\eps}<\infty\}.
\end{equation}

\begin{proposition}\label{prop:L-V}
Let
\begin{equation}
L_t \triangleq \int_0^{t}k(\vP_s)ds+V(\vP_t), \quad t \geq 0.
\end{equation}
Then for every $m \in \mathbb{N}$, $\eps \geq 0$, $\pi \in D$, we
have $L_0 =\E\left\{L_{U_{\eps}\wedge \sigma_m}\right\}$, that is,
\begin{equation}\label{eq:V-V}
V(\vp)=\E\left\{\int_0^{U_{\eps}\wedge
\sigma_m}k(\vP_s)ds+V(\vP_{U_{\eps}\wedge \sigma_m})\right\}.
\end{equation}
\end{proposition}

\begin{proposition}\label{prop:eqp-opti}
The stopping time $U_{\eps}$, which is defined in
(\ref{eq:defn-U-eps}), has bounded $\P$-expectation, for every
$\vp \in D$ and $\eps \geq 0$. More precisely,
\begin{equation}\label{eq:U-eps-bounded}
\E\left\{U_{\eps}\right\}\leq \E\left\{\Theta\right\}+\frac{1}{c},
\quad \vp \in D, \eps \geq 0.
\end{equation}
Moreover, $U_{\eps}$ is $\eps$-optimal for the problem in
(\ref{eq:value-func}); that is
\begin{equation}
\E\left\{\int_0^{U_{\eps}}k(\vP_s)ds+h(\vP_{U_{\eps}})\right\}\leq
V(\vp)+\eps, \quad \vp \in D.
\end{equation}
\end{proposition}

\begin{proof}
Using Proposition~\ref{prop:L-V} and the fact that $V$ is bounded
above by $1$
\begin{equation}
\begin{split}
&1 \geq V(\vp)=\E\left\{\int_0^{U_{\eps}\wedge
\sigma_m}k(\vP_s)ds+V(\vP_{U_{\eps}\wedge \sigma_m})\right\}
\\& \geq \E\left\{\int_0^{U_{\eps}\wedge
\sigma_m}k(\vP_s)ds \right\}= c \E \left\{\left(U_{\eps}\wedge
\sigma_m-\Theta\right)^{+}\right\}\geq  c \E\left\{U_{\eps}\wedge
\sigma_m-\Theta \right\},
\end{split}
\end{equation}
where we used (\ref{eq:peanlty-in-t-Pi}) to derive the second
equality. Applying monotone convergence theorem as $m \uparrow
\infty$, equation (\ref{eq:U-eps-bounded}) follows.

Next, the almost-sure finiteness of $U_{\eps}$ implies
\begin{equation}
V(\vp)=\lim_{m \rightarrow \infty}\E\left\{\int_0^{U_{\eps}\wedge
\sigma_m}k(\vP_s)ds+V(\vP_{U_{\eps}\wedge
\sigma_m})\right\}=\E\left\{\int_0^{U_{\eps}}k(\vP_s)ds+V(\vP_{U_{\eps}})\right\},
\end{equation}
by monotone and bounded convergence theorems, and
Proposition~\ref{prop:L-V}. Since
$V(\vP_{U_{\eps}})-h(\vP_{U_{\eps}}) \geq -\eps$, we have
\begin{equation}
\begin{split}
V(\pi)&=\E\left\{\int_0^{U_{\eps}}k(\vP_s)ds+V(\vP_{U_{\eps}})-h(\vP_{U_{\eps}})+h(\vP_{U_{\eps}})\right\}
\\&\geq
\E\left\{\int_0^{U_{\eps}}k(\vP_s)ds+h(\vP_{U_{\eps}})\right\}-\eps.
\end{split}
\end{equation}
This completes the proof.
\end{proof}

\section{Approximating the Value Function to a Given Level of Accuracy}\label{sec:solution}

In this section, we will describe a numerical procedure that
approximates the value function within any given positive margin,
say $\eps$, and construct $\eps$-optimal stopping strategies.
In the next section, we will give several examples to illustrate
the efficacy of the numerical procedure.

\subsection{Properties of the Stopping Regions}

Let us introduce the stopping and continuation regions for the
problem in (\ref{eq:value-func})
\begin{equation}
\Gamma \triangleq \left\{\vp \in D:
V(\vp)=h(\vp)\right\}, \quad C =D\setminus
\Gamma.
\end{equation}
Taking $\eps=0$ in Proposition~\ref{prop:eqp-opti} implies that
$U_0$ is an optimal stopping time of (\ref{eq:value-func}). From
Remark~\ref{rem:val-func}, we see that an admissible rule to minimize
the Bayes risk in (\ref{eq:value-function}) is to observe the
process $X$ until the process $\vP$ of (\ref{eq:dyn-pi-t}) and
(\ref{eq:dyn-pi-i-t}) enters the stopping region
$\Gamma$.

\begin{remark}\label{rem:con-cont-V}
Since $V$ and $h$ are continuous (the continuity of $V$
follows from Remark~\ref{rem:right-cont-of-paths}),
$\Gamma$ is closed. Moreover, since $V$ is a concave
function (see Corollary~\ref{cor:little-v-n} and
Proposition~\ref{prop:v-V}) and $h$ is linear, $\Gamma$
is a convex set. Indeed, if $\vp_1, \vp_2 \in \Gamma$,
then for any $\alpha \in [0,1]$
\begin{equation}
V(\alpha \vp_1+(1-\alpha) \vp_2) \geq \alpha
V(\vp_1)+(1-\alpha)V(\vp_2)=\alpha h(\vp_1)+
(1-\alpha)h(\vp_2)=h(\alpha \vp_1+(1-\alpha)\vp_2).
\end{equation}
Since $V(\vp) \leq h(\vp)$, for all $\vp \in D$, this equation
implies that $V(\alpha \vp_1+(1-\alpha) \vp_2) =h(\alpha
\vp_1+(1-\alpha)\vp_2)$. Therefore, $\alpha \vp_1+(1-\alpha)\vp_2
\in \Gamma$.
\end{remark}

\begin{proposition}
The stopping region $\Gamma$ is not empty. In particular,
\begin{equation}
\Gamma \supseteqq \left\{\vp=(\pi_1, \cdots, \pi_n, \pi)
\in D: \pi \geq
\frac{\max\{\lambda_0,\lambda_1\}+B}{c+\max\{\lambda_0,\lambda_1\}+B}
\right\},
\end{equation}
in which
\begin{equation}\label{eq:B}
B \triangleq \max_{1 \leq i \leq n} q_{i \Delta}.
\end{equation}

\end{proposition}

\begin{proof}
For $w(\cdot) \ge 0$, using (\ref{eq:distribution-of-first-jump}) and (\ref{eq:sigma-density}),
we write
\begin{align*}
J w(t, \vec{\pi}) &\geq \E\left\{\int_0^{t \wedge
\sigma_1}k(\vec{\Pi}_s)\,ds+1_{\{t<\sigma_1\}}h(\vec{\Pi}_t)\right\}
\geq c \int_0^{t}\pi e^{-\lambda_1 s}\,ds
\\ & \qquad +  c \sum_{i=1}^{n}\pi_i\int_0^{t}
e^{-\lambda_0 s} \left( \int_0^{s}
e^{-(\lambda_1-\lambda_0)(s-u)}f_i(u)du \right)  ds
 + e^{-\lambda_0 t}\sum_{i=1}^{n}\pi_i \int_t^{\infty}f_i(s)\,ds
\\ &\geq
c \int_0^{t}\pi e^{-\lambda_1 s}\,ds+e^{-\lambda_0 t}\sum_{i=1}^{n}\pi_i \int_t^{\infty}f_i(s)\,ds, 
\end{align*}
where $f_i(\cdot)$ is the probability density function of
$\Theta$ given that $M_0=i$, for $i \le n$.

If $\lambda_1 > \lambda_0$, then
\begin{equation}
J w(t, \vec{\pi})  \geq c \int_0^{t}\pi e^{-\lambda_1
s}ds+e^{-\lambda_1 t}\sum_{i=1}^{n}\pi_i
\int_t^{\infty}f_i(s)ds=:K(t,\vp), \quad t \geq 0, \quad \vp \in
D.
\end{equation}
Note that $K(0,\vp)=h(\vp)$. The derivative of $K$ with respect to
$t$
\begin{equation}
\frac{\partial K}{\partial t}(t,\vp)=e^{-\lambda_1 t} \left(c \pi
-\lambda_1 \sum_{i=1}^{n} \pi_i \int_{t}^{\infty}f_i(s)ds-
\sum_{i=1}^{n} \pi_i f_{i}(t) \right).
\end{equation}
Since $f_{i}(t)=dp_{i \Delta}(t)/dt$, Kolmogorov's forward
equation (\ref{eq:kolmogorov}) implies that
\begin{equation}
f_{i}(t) \leq \max_{1 \leq i \leq n}q_{i \Delta}=B.
\end{equation}
Therefore,
\begin{equation}
\frac{\partial K}{\partial t}(t,\vp) \geq 0 \quad \text{if} \quad
\pi \geq \frac{\lambda_1+B}{c+\lambda_1+B}.
\end{equation}
Then for $\pi \geq \frac{\lambda_1+B}{c+\lambda_1+B}$,
\begin{equation}
K(t,\vp) \geq h(\vp) \Rightarrow Jw (t,\vp) \geq h(\vp)
\Rightarrow J_0 w(\vp)=h(\vp).
\end{equation}
Since $V(\cdot) = J_0V(\cdot)$, taking $w=V$ in the last equation, we see that if $\pi \geq
\frac{\lambda_1+B}{c+\lambda_1+B}$, then $\vp=(\pi_1,\cdots,\pi_n,
\pi)\in D$ belongs to $\Gamma$. Similarly, if
$\lambda_0>\lambda_1$ it can be shown that if $\pi \geq
\frac{\lambda_0+B}{c+\lambda_0+B}$, then $\vp=(\pi_1,\cdots,\pi_n,
\pi)\in D$ belongs to $\Gamma$.

\end{proof}

Let us define the optimal stopping and continuation regions for
the problems that we introduced in (\ref{eq:auxiliary}) as
\begin{equation}
\Gamma_m \triangleq \{\vp \in D: V_{m}(\vp)=h(\vp)\}, \quad
\text{and} \quad C_{m}=D\setminus \Gamma_m, \quad m \geq 0.
\end{equation}
Similar arguments as in Remark~\ref{rem:con-cont-V} imply that
$\Gamma_m$ is a closed and convex subset of $D$ for all $m \in
\mathbb{N}_0$. In fact, these sets are ordered, i.e.,
\begin{equation}\label{eq:sets-order}
  \left\{\vp \in D: \pi \geq
\frac{\max\{\lambda_0,\lambda_1\}+B}{c+\max\{\lambda_0,\lambda_1\}+B}
\right\} \subseteq \Gamma \subseteq \cdots \subseteq
\Gamma_{m} \subseteq \cdots \subseteq \Gamma_1 \subseteq \Gamma_0
\equiv D,
\end{equation}
since $V(\cdot) \leq \cdots \leq V_1(\cdot) \leq
V_0(\cdot)=h(\cdot)$.

\subsection{Two Computable $\eps$-Optimal Strategies}

The value function $V(\cdot)$, which is defined in
(\ref{eq:value-function}), can be approximated by the sequence
$\left\{V_m(\cdot) \right\}_{m \in \mathbb{N}_0}$, 
as Proposition~\ref{prop:sequential} suggests. Each element of the
sequence  $\left\{V_m \right\}_{m \in \mathbb{N}_0}$ can be
computed by a successive application of the operator $J_0$, which
is defined in (\ref{eq:defn-J}), to the function $h(\cdot)$, see
(\ref{eq:seq-of-func}) and Proposition~\ref{prop:V-n-epsilon}.
Moreover, the error in approximating $V(\cdot)$ by
$\left\{V_m(\cdot) \right\}_{m \in \mathbb{N}_0}$ can be
controlled. 
Due to Proposition~\ref{prop:sequential}, for
every $\eps>0$, if we choose $\mathcal{M}_{\eps}$ as
\begin{equation}\label{eq:M-eps}
\mathcal{M}_{\eps}=1+\frac{\max\{\lambda_0,\lambda_1\}}{\eps^2}\left(\frac{1}{c}+\E\left\{\theta\right\}\right)
\Rightarrow \|V_{\mathcal{M}}-V\|_{\infty}=\sup_{\vp \in
D}|V_{\mathcal{M}}(\vp)-V(\vp)| \leq \eps, \quad \mathcal{M} \geq
\mathcal{M}_{\eps}.
\end{equation}

In the next section, we will give a numerical algorithm to compute
$V_1, V_2 \cdots$ iteratively. Here, we will describe two
$\eps$-optimal strategies using these functions.

Recall from Proposition~\ref{prop:V-n-epsilon} that
$S_{m}^{\eps}$, $m \geq 1$ are $\eps$-optimal stopping times for
the problem in (\ref{eq:auxiliary}). For a fixed $\eps>0$, if
we choose $\mathcal{M} \geq \mathcal{M}_{\eps/2}$, we have 
$\|V_{\mathcal{M}}-V\|_{\infty} \leq \eps/2$. Then
$S_{\mathcal{M}}^{\eps/2}$ is $\eps-$optimal for $V(\cdot)$
since
\begin{equation}
\E\left\{\int_0^{S_{\mathcal{M}}^{\eps/2}}k(\vec{\Pi}_t)dt+h(\vec{\Pi}_{S_{\mathcal{M}}^{\eps/2}})\right\}
\leq V_{\mathcal{M}}(\vec{\pi})+\frac{\eps}{2} \leq V(\vp)+\eps,
\quad \vp \in D.
\end{equation}
Note that $S_{\mathcal{M}}^{\eps/2}$ is not a hitting time. In (\ref{eq:defn-S-eps}), it prescribes to wait until the minimum of
$r_{\mathcal{M}-1}^{\eps/4}(\vp)$ and the first jump time
$\sigma_1$ of the process $X$. If
$r_{\mathcal{M}-1}^{\eps/4}(\vp)<\sigma_1$, then we stop.
Otherwise, the probabilities are updated to $\vP_{\sigma_1}$ and
we wait until the minimum of
$r_{\mathcal{M}-1}^{\eps/4}(\vP_{\sigma_1})$ and the next jump
time $\sigma_2=\sigma_1 \circ \theta_{\sigma_1}$ of the process
$X$. If $r_{\mathcal{M}-1}^{\eps/4}(\vP_{\sigma_1})$ comes first
we stop. Otherwise we continue as before. We finally stop at
the $M$th jump time if not before.

We can also give an $\eps$-optimal strategy that is a hitting
time. Let us define
\begin{equation}\label{eq:U-M-eps}
U_{\eps/2}^{(\mathcal{M})} \triangleq \inf\left\{t \geq 0:
h(\vP_t) \leq V_{\mathcal{M}}(\vP_t)+\frac{\eps}{2} \right\}.
\end{equation}
Following the same arguments as in the proof of
Proposition~\ref{prop:eqp-opti} this stopping time can be shown to
be an $\eps/2$-optimal stopping time for $V_{\mathcal{M}}(\cdot)$,
which in turn implies that it is an $\eps-$optimal stopping time
for $V(\cdot)$.

\subsection{An Algorithm Approximating the Value Function}
\label{sec:algorithm}
Note that if the hitting time of $t \rightarrow x_{\Delta}(t,\vp)$ to the region $\Gamma$
is uniformly bounded by some $t^* < \infty$, then the minimization problem in computing
$V_{m+1}(\vec{\pi}) 
=\inf_{t \in [0,\infty]}JV_m(t,\vp)$ can be restricted to the
compact interval $ [0,t^*]$ thanks to Corollary~\ref{cor:r-m}.
Remark~\ref{rem:time-bound} constructs a uniform bound $t^*$ when
the parameters of the problem satisfy
$\tilde{B}-\lambda_1+\lambda_0 \geq 0$, in which $\tilde{B}$
defined as
\begin{equation}\label{eq:Btilde}
\tilde{B} \triangleq \min_{1 \leq i \leq n}q_{i \Delta}.
\end{equation}


\begin{remark}\label{rem:time-bound}
The hazard rate of the prior distribution
of $\theta$ satisfies
$\eta(t)\geq \tilde{B}$ (see (\ref{eq:alt-for-h})).
Moreover, from (\ref{eq:pi-h}), we have 
\begin{equation}
\frac{d x_{\Delta}(t,\vp)}{dt}=(\eta(t)-(\lambda_1-\lambda_0)x_{\Delta}(t,\vp) )(1-x_{\Delta}(t,\vp)),
\; \quad x_{\Delta}(0, \vp)=\pi,
\end{equation}
where $x_{\Delta}$ is defined in (\ref{eq:dyn-x}). Let $\tilde{x}(t)$ be the solution of the differential equation
\begin{equation}\label{eq:auxi-diff}
\frac{d \tilde{x}(t)}{dt}=(\tilde{B}-(\lambda_1-\lambda_0)\tilde{x}(t))(1-\tilde{x}(t)),
\qquad \text{with }\; \tilde{x}(0)=0.
\end{equation}
A simple comparison argument shows that $x_{\Delta}(t) \geq x(t)$,
for all $t \geq 0$ when $\tilde{B}-\lambda_1+\lambda_0 \geq 0$.
The solution to (\ref{eq:auxi-diff}) can be written as
\begin{equation}\label{eq:explicit-x}
x\tilde{}(t)= \begin{cases}
\frac{\frac{\tilde{B}}{\tilde{B}-\lambda_1-\lambda_0}\left(1-\exp\left((\tilde{B}-\lambda_1+\lambda_0)t\right)\right)}
{1+\frac{\tilde{B}}{\tilde{B}-\lambda_1-\lambda_0}\left(1-\exp\left((\tilde{B}-\lambda_1+\lambda_0)t\right)\right)}
& \text{if \,\, $\tilde{B}-\lambda_1+\lambda_0 \ne 0$},
\\ \frac{\tilde{B}t}{1+\tilde{B}t} & \text{if \,\,
$\tilde{B}-\lambda_1+\lambda_0=0$}.
\end{cases}
\end{equation}
When $\tilde{B}-\lambda_1+\lambda_0 \geq 0$, let us denote
\begin{equation}
\hat{x} \triangleq
\frac{\max\{\lambda_0,\lambda_1\}+B}{c+\max\{\lambda_0,\lambda_1\}+B},
\end{equation}
where $B$ is given in (\ref{eq:B}). Let $t^*(\vp)  \triangleq \inf \{ t > 0;  x_{\Delta}(t,\vp)=\hat{x} \}$,
then using
(\ref{eq:explicit-x}), it can be easily verified that
\begin{equation}
\label{eq:uniform-bound-on-hitting-time}
t^*(\vp) \leq t^* \triangleq \begin{cases}
\frac{1}{\tilde{B}-\lambda_1-\lambda_0} \log\left(\frac{\tilde{B}+
\hat{x} (\lambda_0-\lambda_1)}{\tilde{B}(1-\hat{x})}\right), & \text{if
\,\, $\tilde{B}-\lambda_1+\lambda_0>0$,}
\\ \frac{\hat{x}}{1-\hat{x}}\frac{1}{\tilde{B}} , & \text{if
\,\, $\tilde{B}-\lambda_1+\lambda_0=0$},
\end{cases}
\end{equation}
for all $\vp \in D$.
\end{remark}

\begin{remark}
If $\tilde{B}-\lambda_1+\lambda_0 < 0$, then $t^{*}(\vp)$ defined
in Remark~\ref{rem:time-bound} may be $\infty$ for some $\vp \in D$.
\end{remark}

When $\tilde{B}-\lambda_1+\lambda_0 < 0$, it is still possible to
restrict the minimization problem in $V_{n+1}(\vp) = \inf_{t \in [0,\infty]} J V_n
(t,\vp)$ to a compact interval and control the error arising from
this. Note that for any $w(\cdot) \leq 1$, we have
\begin{equation}
\begin{split}
& \sup_{\vp \in D}|Jw (t,\vp)-J w (\infty, \vp)| \\&\leq c
\int_t^{\infty}\P\{s \leq \sigma_1\}
ds+\int_t^{\infty}\P\{\sigma_1 \in ds\}
Sw(\vec{x}(s,\vec{\pi}))+h(\vec{x}(t,\pi))\P\{t<\sigma_1\}
\\& \leq c
\int_t^{\infty}\P\{s \leq \sigma_1\} ds+ 2 \P\left\{\sigma_1 \geq
t \right\} \leq c \int_t^{\infty}e^{-\lambda_0 s}ds+2
e^{-\lambda_0 t} \leq
\left(\frac{c}{\lambda_0}+2\right)e^{-\lambda_0 t},
\end{split}
\end{equation}
where the first inequality follows from (\ref{eq:exp-for-J-w}),
and the second one follows from the fact that $V_m \leq h \leq 1$,
the third one follows from $\P\{\sigma_1 \geq t\} \leq
e^{-\lambda_0 t}$, which is a direct consequence of
(\ref{eq:distribution-of-first-jump}). Then, denoting
\begin{equation}\label{eq:t-delta}
t(\delta)\triangleq -\frac{1}{\lambda_0} \log
\left(\frac{\delta}{4+2c/\lambda_0}\right),
\end{equation}
we obtain
\begin{equation}\label{JV-t1-t2}
|Jw (t_1,\vp)-Jw(t_2,\vp)|\leq
|Jw(t_1,\vp)-Jw(\infty,\vp)|+|Jw(t_2,\vp)-Jw(\infty,\vp)|\leq
\delta.
\end{equation}
for any $t_1,t_2 \geq t(\delta)$. Letting
\begin{equation}\label{eq:J-0-delta}
J_{0,t}w(\vp) \triangleq \inf_{s \in [0,t]}J w(s,\vp), \quad
\text{for every bounded $w:D\rightarrow \R$, and $t \geq 0$, $\vp \in
D$},
\end{equation}
we get
$ \sup_{\vp \in D}|J_{0,t(\delta)}w(\vp)-J_0 w(\vp)| \leq \delta$.
Now, let us define a new sequence of functions as
\begin{equation}
V_{\delta,0}(\vp) \triangleq h(\vp) \quad \text{and} \quad
V_{\delta,m+1}(\vp) \triangleq J_{0, t(\delta)}V_{\delta,m}(\vp),
\quad \vp \in D.
\end{equation}
\begin{proposition}\label{prop:new-approx}
For every $\delta>0$, $m \geq 0$, we have
\begin{equation}\label{eq:new-bounds}
V_m(\vp) \leq V_{\delta,m}(\vp) \leq m \delta+V_m (\vp), \qquad \vp
\in D.
\end{equation}
\end{proposition}
\begin{proof}
For $m=0$ we have $V_{\delta,0}(\cdot)=V_0(\cdot)=h(\cdot)$, 
by construction. Now, suppose that (\ref{eq:new-bounds}) holds
for some $m \geq 0$. Then
\begin{equation}
V_{m+1}(\vp)=J_0 V_{m}(\vp) \leq J_0 V_{\delta,m}(\vp) \leq J_{0,
t(\delta)}V_{\delta,m}(\vp) \leq V_{\delta,m+1} \qquad \vp \in D,
\end{equation}
which proves the first inequality in (\ref{eq:new-bounds}) when we
substitute $m$ with $m+1$. The first inequality follows from the
induction hypothesis and Remark~\ref{rem:monotone}. The second
inequality follows from (\ref{eq:J-0-delta}). 

Let us now prove the second inequality in (\ref{eq:new-bounds}) when $m$ is replaced by
$m+1$. Observe that $V_{\delta,m}(\vp) \leq h(\vp)$, $\vp \in D$.
Then
\begin{equation}
\begin{split}
V_{\delta, m+1}(\vp)&=\inf_{t \in [0,t(\delta)]}J
V_{\delta,m}(t,\vp)\leq \inf_{t \in [0,\infty]}J
V_{\delta,m}(t,\vp)+\delta \\ &\leq \inf_{t \in [0,\infty]}\left[J V_m
(t,\vp)+m \delta \int_0^{t}\P\{\sigma_1 \in ds\}\right]+\delta
\\& \leq V_{m+1}(\vp)+ m \delta \int_0^{\infty}\P\{\sigma_1 \in
ds\}+\delta \leq  V_{m+1}(\vp)+(m+1)\delta,
\end{split}
\end{equation}
where the second inequality follows from the induction hypothesis
and the definition of the operator $J$.
\end{proof}

When $\tilde{B}-\lambda_1+\lambda_0 < 0$, using
Proposition~\ref{prop:new-approx} we can approximate the value function $V(\cdot)$
with the functions
$\left\{V_{\delta,m}(\cdot)\right\}_{\delta>0,m\geq 1}$. There is
an extra error, because we truncate at $t(\delta)$, but this can
be compensated by increasing the number of iterations. Let us
define
\begin{equation}\label{eq:tilde-M-eps}
\tilde{\mathcal{M}}_{\eps}\triangleq
1+\frac{1}{\eps^2}\left[1+\sqrt{\left(\frac{1}{c}+\E\{\Theta\}\right)\max\{\lambda_0,\lambda_1\}}\right],
\quad \text{and} \quad \delta_{\eps} \triangleq
\frac{1}{\tilde{\mathcal{M}}_{\eps}\sqrt{\tilde{\mathcal{M}}_{\eps}-1}}.
\end{equation}
Then for every $\mathcal{M} \geq \tilde{\mathcal{M}}_{\eps}$ and
$\delta \leq \delta_{\eps}$ we have
\begin{equation}
\|V_{\delta,\mathcal{M}}-V\|_{\infty} \leq
\|V_{\delta,\mathcal{M}}-V_{\mathcal{M}}\|+\|V_{\mathcal{M}}-V\|
\leq \mathcal{M} \delta+
\sqrt{\left(\frac{1}{c}+\E\{\Theta\}\right)\frac{\max\{\lambda_0,
\lambda_1\}}{\mathcal{M}-1}} \leq \eps,
\end{equation}
where we used Propositions~\ref{prop:sequential} and
\ref{prop:new-approx}. In other words, by applying the
operator $J_{0,t(\delta_{\eps})}$ to the function $h(\cdot)$
$\tilde{M}_{\eps}$ times, we obtain an approximation of $V(\cdot)$
within $\eps$-closeness on $D$.

Similar arguments as in Section~\ref{sec:sequential} can be repeated to show that each $V_{\delta,m}$, for $m \ge 1$, is continuous and concave on $D$. Moreover, we can still define $\eps$-optimal rules using the function $V_{\delta,\mathcal{M}}$.
Particularly, let us define the stopping time
\begin{equation}
U_{\eps/2}^{(\mathcal{M},\delta)} \triangleq \inf\left\{t \geq 0:
h(\vP_t) \leq V_{\delta,M}(\vP_t)+\frac{\eps}{2} \right\}.
\end{equation}
When we take $\mathcal{M}=\tilde{\mathcal{M}}_{\eps/2}$ and
$\delta=\delta_{\eps/2}$, this stopping time becomes an
$\eps-$optimal stopping time for the problem in
(\ref{eq:value-function}). This follows using the same arguments
as in the proof of Proposition~\ref{prop:eqp-opti}.

Finally, we conclude this section with the following numerical algorithm
summarizing the results presented here in order to approximate $V(\cdot)$.

\textbf{Algorithm.}

1) If $\tilde{B}-\lambda_1+\lambda_0 \geq 0$, then choose
$\mathcal{M} \geq \mathcal{M}_{\eps}$. $\tilde{B}$ is given by
(\ref{eq:Btilde}) and $M_{\eps}$ is given by (\ref{eq:M-eps}).

1') On the other hand if $\tilde{B}-\lambda_1+\lambda_0 < 0$, then
choose $\mathcal{M} \geq \tilde{\mathcal{M}}_{\eps}$ and $\delta
\leq \delta_{\eps}$, in which $\tilde{\mathcal{M}}_{\eps}$ and
$\delta_{\eps}$ are given by (\ref{eq:tilde-M-eps}).

2) Set $V_0(\cdot)=h(\cdot)$.

2') Set $V_{\delta,0}(\cdot)=h(\cdot)$.

3) Calculate $V_{m+1}(\vp)=\min_{t \in [0,t^{*}(\vp)]}J
V_m(t,\vp)$, $\vp \in D$, in which $t^{*}(\vp)$ is
given in Remark~\ref{rem:time-bound} (see also (\ref{eq:uniform-bound-on-hitting-time})).
%

3') Calculate $V_{\delta,m+1}=\inf_{t \in [0,t(\delta)]}J
V_{\delta,n}(t,\vp)$, $\vp \in D$, in which $t(\delta)$ is defined
in (\ref{eq:t-delta}).

4) Repeat step 3 until $m=\mathcal{M}+1$.

4') Repeat step 3' until $m=\mathcal{M}+1$.

If  $\tilde{B}-\lambda_1+\lambda_0 \geq 0$, our algorithm returns
$V_\mathcal{M}$, which satisfies $\|V_\mathcal{M}-V\| \leq \eps$.
On the other hand if $\tilde{B}-\lambda_1+\lambda_0< 0$, the
algorithm returns $V_{\delta,\mathcal{M}}$, which satisfies
$\|V_{\delta,\mathcal{M}}-V\|\leq \eps$.

\section{Examples}\label{sec:examples}
In this section, we provide examples illustrating the use of the numerical algorithm presented above for negligible $\varepsilon$-values.
\subsection{Mixed Erlang distribution}
In (\ref{eq:matrix}) let us take a particular form for $\mathcal{A}$ where all entries are zero except  $q_{ii} = - \lambda$, $ q_{i,i+1} = \lambda$ for
some rate $\lambda >0 $, and for $i=1,\ldots,n$. Then, starting
from any non-absorbing state $i$, the process $M$ visits all the states $i+1, i+2, \ldots$
until it eventually hits the absorbing state $\Delta$.  In other words, conditioned on any initial non-absorbing state $i$,
the disorder time has Erlang distribution with the shape index $n-i+1$ and rate $\lambda$. In this case the distribution of $\Theta$ can be explicitly given as
\begin{align*}
F_{\vp}(t) = \P\{ \Theta \le t \} = \pi + \sum_{i \ne \Delta}
\pi_i \cdot  \int_0^t f_i(u) du, \qquad \text{in terms of} \qquad
f_i (t) \triangleq \frac{ \lambda^{n+1-i} t^{n-i} }{ (n-i)! }
  e^{-\lambda t} ,
\end{align*}
for $i \le n$. Moreover, the components of the deterministic path $\vec{x}(\cdot, \cdot)$ have the explicit forms
\begin{align*}
x_i(t,\vp) = \frac{\sum_{j=1}^i  \pi_j  e^{-\lambda t} \frac{
(\lambda t)^{i-j} }{ (i-j)! }  }{ \left( \sum_{k=1}^{n}
\sum_{j=1}^{k} \pi_j  e^{-\lambda t} \frac{  (\lambda t)^{k-j}}{
(k-j)! } \right) + e^{-(\lambda_1 - \lambda_0)t} \left( \pi +
\sum_{k=1}^{n} \pi_k \int_0^t e^{(\lambda_1 - \lambda_0)u} f_k (u)
du \right) },
\end{align*}
for $i \le n$, and for the $n+1$'st component $x_{\Delta}$ we have
\begin{align*}
x_{\Delta}(t,\vp) = \frac{e^{-(\lambda_1 - \lambda_0)t} \left( \pi
+ \sum_{k=1}^{n} \pi_k \int_0^t e^{(\lambda_1 - \lambda_0)u} f_k
(u) du   \right)   }{ \left( \sum_{k=1}^{n} \sum_{j=1}^{k} \pi_j
e^{-\lambda t}  \frac{ (\lambda t)^{k-j} }{ (k-j)! } \right) +
e^{-(\lambda_1 - \lambda_0)t} \left( \pi + \sum_{k=1}^{n} \pi_k
\int_0^t e^{(\lambda_1 - \lambda_0)u} f_k (u) du   \right) }.
\end{align*}
Using these expressions (and assuming $\pi \ne 1$), it can be shown that if $\lambda - \lambda_1+ \lambda_0 \ge 0 $ then $x_{\Delta}(t,\vp) \to 1$ as $t \to \infty$. Otherwise, we have
\begin{align}
\label{convergence-points-in-mixed-erlang}
\lim_{t \to \infty} x_n(t,\vp) = \frac{\lambda_1 - \lambda_0 - \lambda}{\lambda_1 - \lambda_0 } \, , \quad \text{ and} \quad \lim_{t \to \infty} x_{\Delta}(t,\vp)  = \frac{  \lambda }{\lambda_1 - \lambda_0 } .
\end{align}
Due to the explicit form of the paths $t \mapsto \vx(\cdot, \cdot)$, the steps described in the numerical algorithm above are easier to carry. The Figure 1 below illustrates examples on two different problems where there are two transient states. 

In Panels (a) and (b) of Figure 1, we see the sample path behavior and the value function of a problem where the parameters are $\lambda_0= 6$, $\lambda_1=5 $, $\lambda=3$, $c=1$. Panel (a) presents the behavior of the paths $t \to \vec{x} (t, \vp) $ for a number of different starting points. We also plot a sample path of $\vP_t$ starting from a particular point. Since $\lambda_0 > \lambda_1$, the $n+1$'st component $x_{\Delta}$ of $\vx$ is increasing. In other words, as long as we do not observe any arrival, we tend to assign more likelihood the event that the disorder has happened by then.  On the other hand, when we observe an arrival, we decrease this likelihood. Moreover, since $\lambda - \lambda_1 + \lambda_0  \geq 0$ we see that the paths of $\vx$ converge asymptotically to the point $(0,0,1)$ as indicated above. In this case, we use steps (1), (2), (3) and (4) of the algorithm that is presented at end of Section~\ref{sec:algorithm} to approximate the value function to a given order of accuracy. Thanks to the properties of the approximating sequence (see Section~\ref{sec:sequential}), properties such as concavity of the value function and the convexity of the optimal stopping boundary are  preserved by our approximation.
Panel (b), on the right, illustrates the (approximated) value function defined on the state space $D$ of $\vP$.  As the figure shows, the value function $V(\cdot)$ is non-negative and concave on $D$, and there exists a region on the neighborhood of the point $(0,0,1)$ where it coincides with the terminal reward function $h(\cdot)$.  As indicated in Section 4, an ($\varepsilon$-)optimal strategy then implies that one observes the counting process $X$, and update the process $\vP$ continuously until $\vP$ enters the region $\Gamma$. At this time, we stop and declare that the disorder has happened by then.

\begin{center}
\begin{tabular*}{\textwidth}
    {@{\extracolsep{\fill}}cc}
    \includegraphics[scale=0.4]{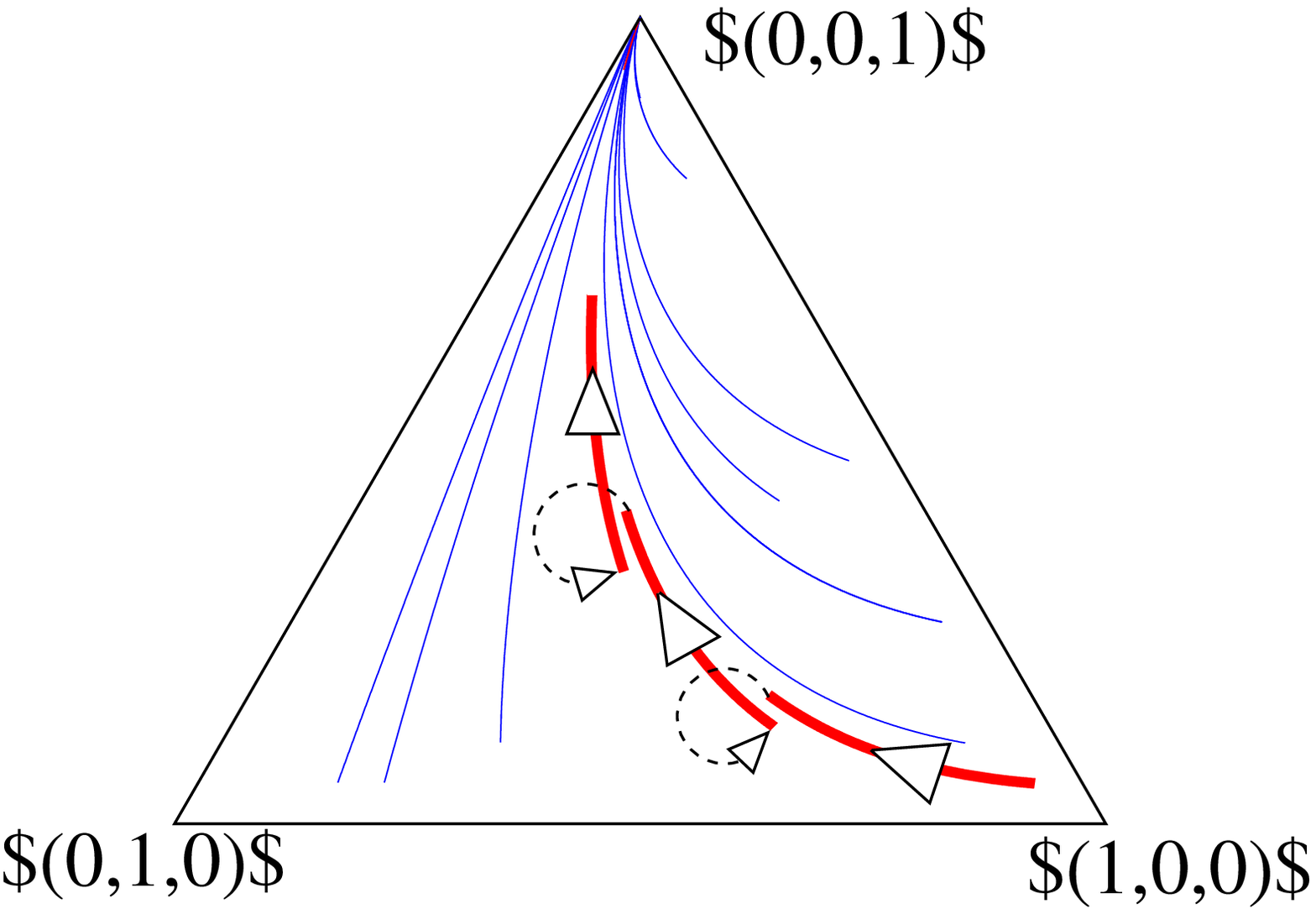}&
     \includegraphics[scale=0.65]{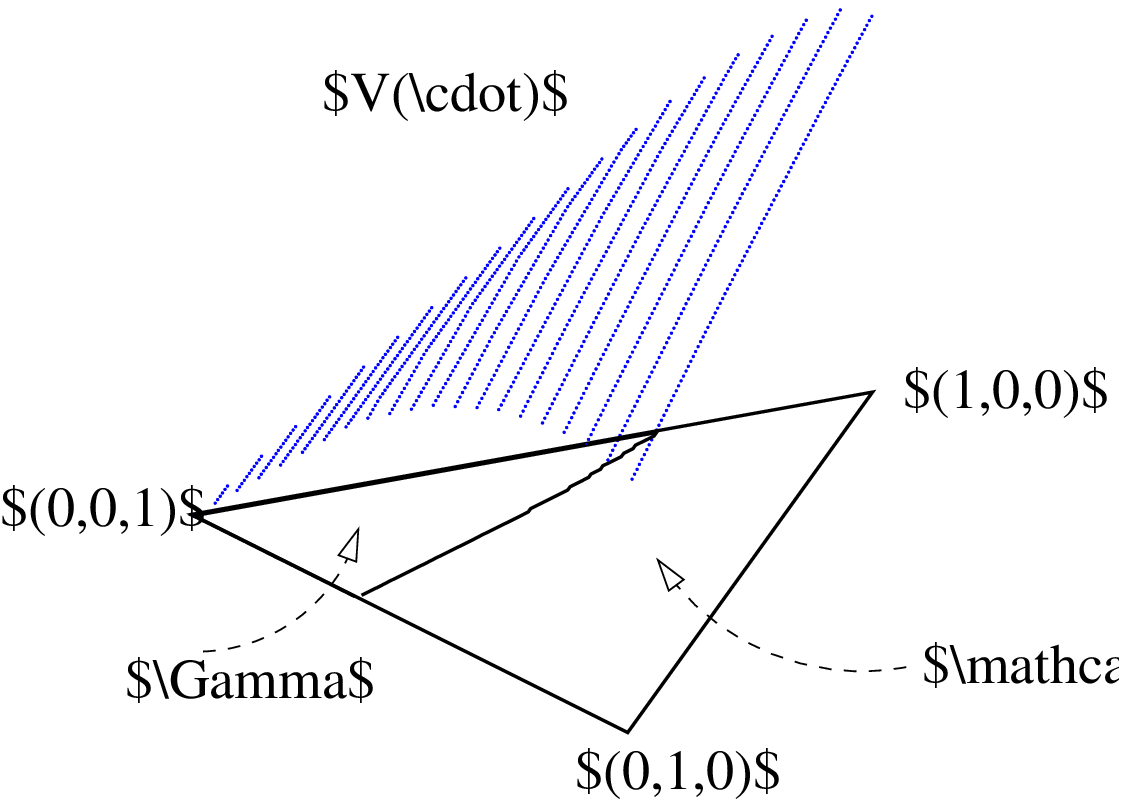}  \\
    (a) & (b) \\ \text{ } & \text{ } \\
   \includegraphics[scale=0.4]{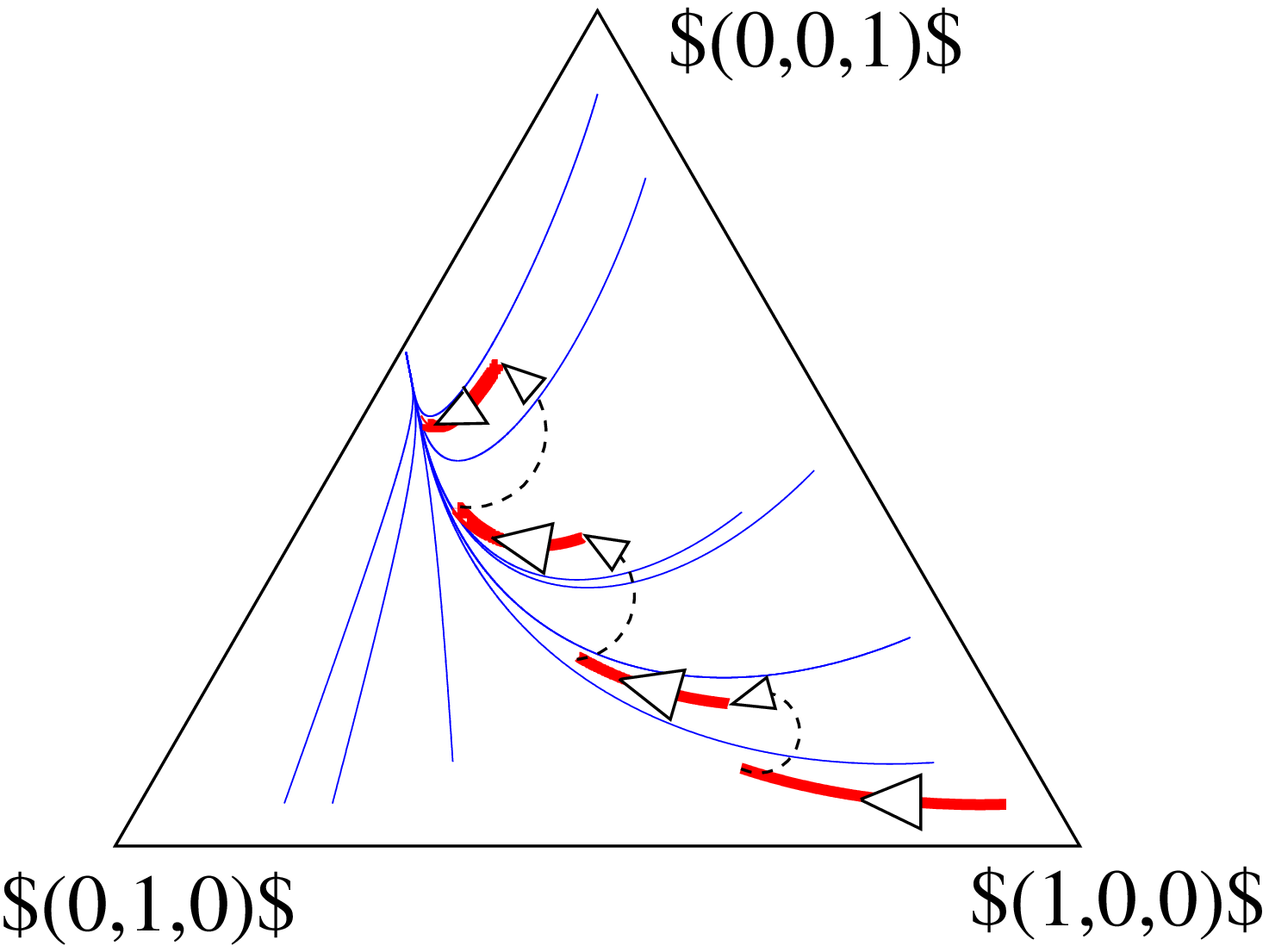}&
    \includegraphics[scale=0.65]{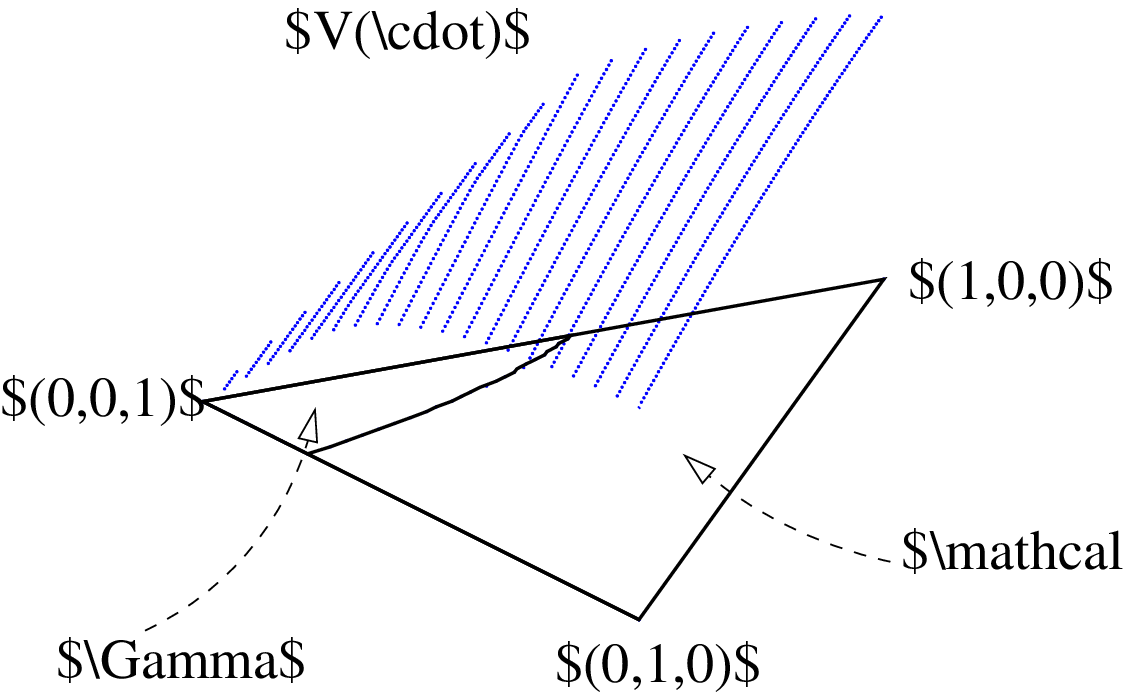}  \\ 
    (c) & (d) \\ \text{ } & \text{ } \\
\end{tabular*}
\emph{\textbf{Figure 1}: Examples with mixed Erlang prior distributions. Panels (a) and (b) correspond to a problem with $\lambda_0= 6$, $\lambda_1=5 $, $\lambda=3$, $c=1$. Panel (a) represents the sample path behavior of $t \mapsto \vx(t,\vp)$ and $t \mapsto \vP_t$. The continuous curves are possible sample paths of $\vx$ starting from different points. The discontinuous path with arrows indicates the behavior of $\vP$. As indicated in Section 2, between two jumps, the process $\vP$ follows the deterministic curves of $\vx$, and at jump times it switches from one curve to another. Panel (b) gives the value function $V(\cdot)$ and the stopping ($\mathcal{C}$) and continuation ($\Gamma$) regions. Similarly, Panels (c) and (d) correspond to another problem where $\lambda_0= 5$, $\lambda_1=10 $, $\lambda=3$, $c=1$.}
\end{center}

Panels (c) and (d), on the other hand, correspond to another sample problem where $\lambda_0= 5$, $\lambda_1=10 $, $\lambda=3$, $c=1$. In this case, we have $\lambda_1 > \lambda + \lambda_0$. Therefore, the paths $t \mapsto \vx(t,\vp)$ are asymptotically converging to the point $(0,0.4,0.6)$ as indicated in (\ref{convergence-points-in-mixed-erlang}). Moreover, since $\lambda_1 >  \lambda_0 $, the $n+1$'st component $\Pi_t$ of $\vP$ moves closer to the point $(0,0,1)$ at jump times. In this case we use steps (1'), (2'), (3') and (4') of the algorithm that we presented at the end of Section~\ref{sec:algorithm} to approximate the value function. In Panel (d), we verify our claims in Section~\ref{sec:sequential} again. That is; the value function is a concave function and the stopping region $\Gamma$ is a convex region around the point $(0,0,1)$.

\subsection{Hyperexponential distribution}
Let us here reconsider the formulation of the classical Poisson disorder problem with exponential prior distribution, and let us assume that the rate of this exponential distribution is not known precisely. Rather there are $n$ possible rates $\mu_1, \mu_2, \ldots, \mu_n$ with prior likelihoods $(\pi_1, \ldots , \pi_n, \pi)$, and the aim is to detect the change time $\Theta$ by minimizing $R_{\tau}(\vec{\pi})$ in (\ref{eq:bayes-risk}).

This problem can be modeled as a special case of phase-type Poisson disorder problem if we take column vector $r$ in (\ref{eq:matrix}) in the form $r =  [\mu_1, \mu_2, \ldots, \mu_n]'$ for $\mu_i >0$ for $i=1,\ldots, n$. Moreover we let the matrix $R$ in (\ref{eq:matrix}) be $R= - r' \cdot I$, where $I$ is $n \text{x} n$ identity matrix.
In this case, if the process $M$ starts from a transient state, it is absorbed to the state $\Delta$ at the first transition time, and conditioned on the initial state $i$ the hitting time has exponential distribution with parameter $\mu_i$.

In this case, by direct computation it can be shown that the deterministic paths $\vx(\cdot, \cdot)$ has the form
\begin{align}
\label{deterministic-paths-for-hyperexp-case} x_i(t,\vp) = \frac{
\pi_i  e^{-\mu_i t}     }{ \left( \sum_{k=1}^{n}  \pi_k  e^{-\mu_k
t}  \right) +e^{-(\lambda_1 - \lambda_0)t} \left( \pi +
\sum_{k=1}^{n} \pi_k \int_0^t e^{(\lambda_1 - \lambda_0)u} f_k (u)
du \right)  },
\end{align}
for $1 \le i \le n$ and
\begin{align*}
x_{\Delta}(t,\vp) =\frac{ e^{-(\lambda_1 - \lambda_0)t} \left( \pi
+ \sum_{k=1}^{n} \pi_k \int_0^t e^{(\lambda_1 - \lambda_0)u} f_k
(u) du   \right)   }{ \left( \sum_{k=1}^{n}  \pi_k  e^{-\mu_k t}
\right) + e^{-(\lambda_1 - \lambda_0)t} \left( \pi +
\sum_{k=1}^{n} \pi_k \int_0^t e^{(\lambda_1 - \lambda_0)u} f_k (u)
du \right) },
\end{align*}
where $f_k (u) = \mu_k e^{- \mu_k u}$, for $1 \le k \le n$.

Without loss of generality let us assume that $\mu_1> \mu_2 > \ldots > \mu_n$. Then, on $\{ \vp \in D: \pi_n \ne 0\}$, the path $x_i(t,\vp)$ goes to $0$ as $t \to \infty$ for $i=1,\ldots, n-1$. If $\mu_n - \lambda_1 + \lambda_0 \ge  0$, then $x_i(t,\vp)$ converges to $1$ asymptotically, otherwise we have
\begin{align}
\label{convergence-points-in-hyperexponential}
\lim_{t \to \infty} x_n(t,\vp) = \frac{ \lambda_1 -\mu_n - \lambda_0  }{\lambda_1 - \lambda_0 } , \quad \text{ and} \quad \lim_{t \to \infty} x_{\Delta}(t,\vp)  = \frac{  \mu_n }{\lambda_1 - \lambda_0 }
\end{align}
On the other hand, on the region $\{ \vp \in D: \pi_n = 0, \pi_{n-1} \ne 0\}$, the above statements hold by replacing $n$ and $\mu_n$ with $n-1$ and $\mu_{n-1}$ respectively, and so on. If a non-absorbing state has the initial likelihood $\pi_i=0$, then $\Pi^{(i)}_t =0$, for all $t \ge 0$ by (\ref{deterministic-paths-for-hyperexp-case}) and (\ref{eq:rel-pi-x}). Indeed, since the disorder occurs at the first transition time, this state can be eliminated from the problem.

Figure 2 presents two numerical examples with two transient states. In Panels (a) and (b), we see the value function and the paths $\vx(\cdot, \cdot)$ of a problem where the parameters are $\mu_1=3$, $\mu_2=2$, $\lambda_0=2$, $\lambda_1=1$, $c=1.5$. Between two jumps, the process $\vP$ follows the paths $t \mapsto \vx(t, \vp)$, which are converging to the point $(0,0,1)$ asymptotically. Moreover since $\lambda_1 < \lambda_0$, the process $\vP$ jumps away from this point, and we decrease the conditional likelihood of the disorder event at arrival times of $X$. In this case, we use steps (1), (2), (3) and (4) of the algorithm that is presented at end of Section~\ref{sec:algorithm} to approximate the value function. In Panel (b), we observe that the value function is concave, and the stopping region is a convex region with non-empty interior around the point $(0,0,1)$ as indicated in Section~\ref{sec:sequential}.

\begin{center}
\begin{tabular*}{\textwidth}
    {@{\extracolsep{\fill}}cc}
    \includegraphics[scale=0.4]{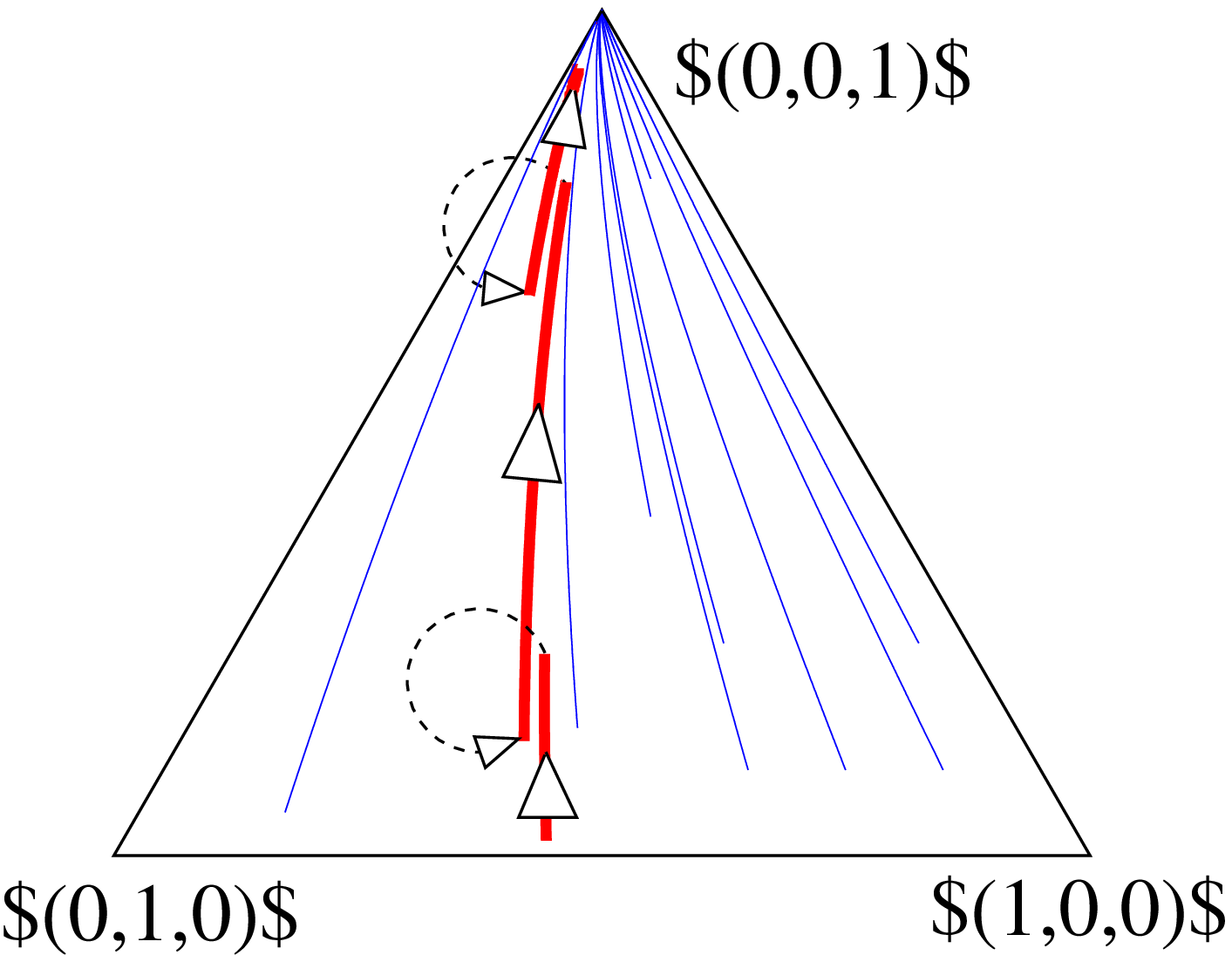}&
     \includegraphics[scale=0.65]{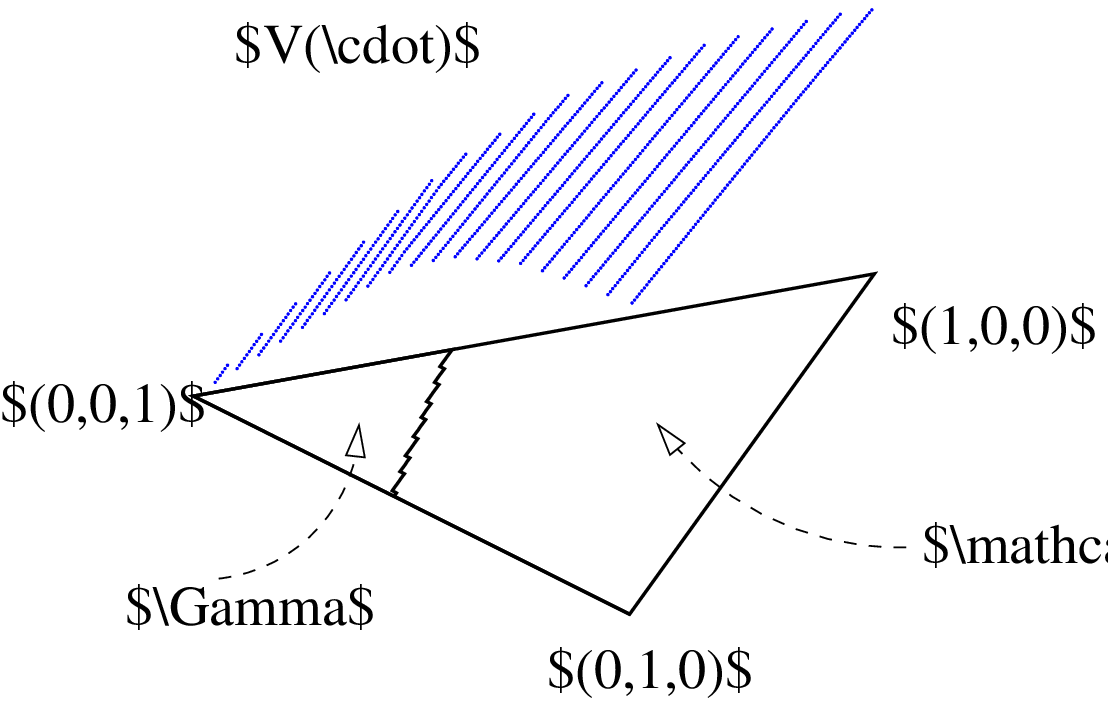} \\
    (a) & (b) \\ \text{ } & \text{ } \\
    \includegraphics[scale=0.4]{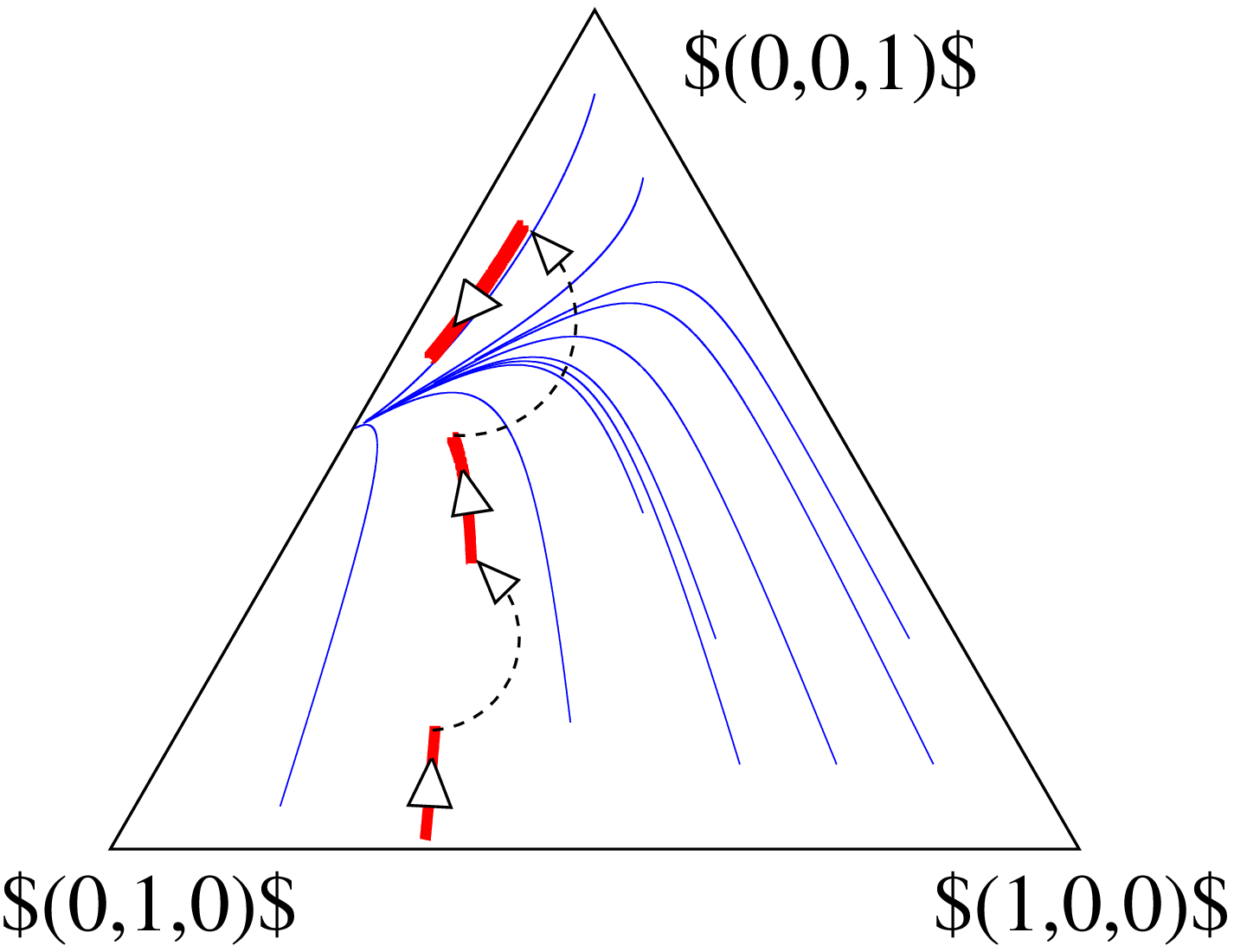}&
     \includegraphics[scale=0.65]{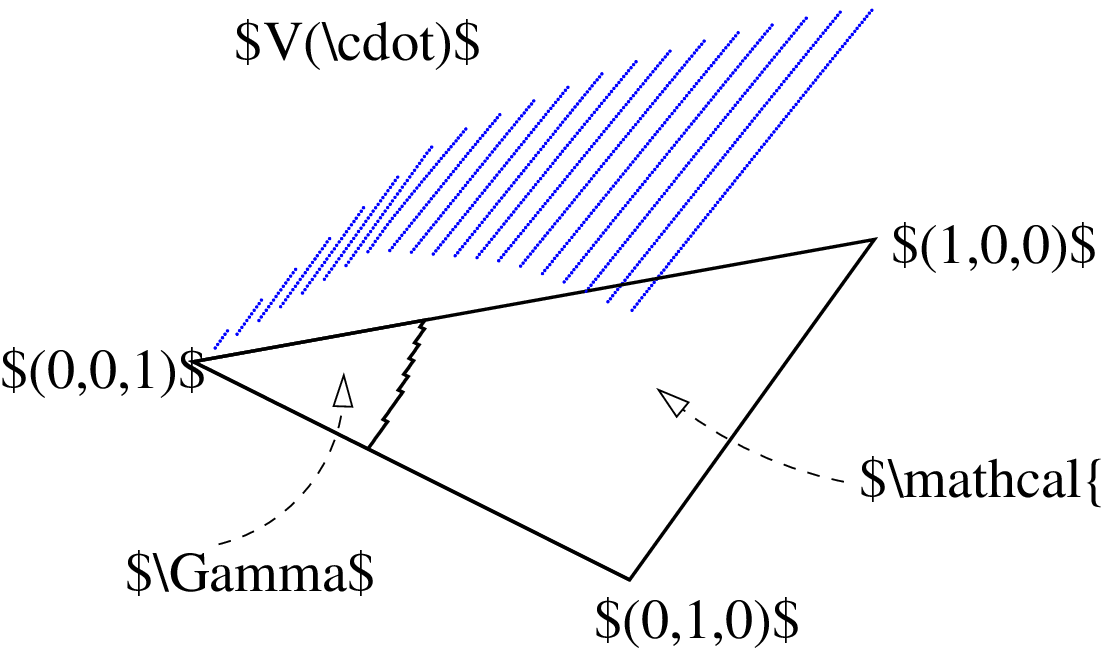} \\
    (c) & (d) \\ \text{ } & \text{ } \\
\end{tabular*}
\emph{\textbf{Figure 2}: Examples with hyper-geometric prior distributions. In panels (a) and (b) we see the sample path properties and the value function of a problem 
where $\mu_1=3$, $\mu_2=2$, $\lambda_0=2$, $\lambda_1=1$, $c=1.5$. Continuous paths in Panel (a) are the paths of $t \mapsto \vx(t,\vp)$ for different starting points.  The discontinuous path with arrows is a sample path of $t \mapsto \vP_t$. Panel (b) illustrates the value function $V(\cdot)$ and the stopping ($\mathcal{C}$) and continuation ($\Gamma$) regions. Similarly, in Panels (c) and (d) we see another problem with $\mu_1=3$, $\mu_2=2$, $\lambda_0=2$, $\lambda_1=6$, $c=1.5$.}
\end{center}

In panels (c) and (d) we have another problem whose parameters are $\mu_1=3$, $\mu_2=2$, $\lambda_0=2$, $\lambda_1=6$, $c=1.5$. In this case we have $\mu_n - \lambda_1 + \lambda_0 < 0$. Hence, in accordance with (\ref{convergence-points-in-mixed-erlang}) we see that the paths $t \mapsto \vx(t,\vp)$ converge to the point $(0, 0.5, 0.5)$. Also, since $\lambda_1 > \lambda_0$, at the jump times of $X$, the process $\vP$ jumps towards the point $(0,0,1)$ and the conditional probability of the disorder event is increased.  In this case we use steps (1'), (2'), (3') and (4') of the algorithm that we presented at the end of Section~\ref{sec:algorithm} to approximate the value function. In Panel (d), we verify once again the concavity of the value function and convexity of the stopping region around the point $(0,0,1)$.

Non-smooth behavior of the value function on the region where $\pi_1 =0$ is in accordance with Lemma 7.1 of \cite{ds}. On this region, the problem is essentially with one non-absorbing state. The point $\vp= (\,0,1- \mu_2 / (\lambda_1 - \lambda_0) , \mu_2 / (\lambda_1 - \lambda_0) \,) = (0, 0.5, 0.5)$ falls into the continuation, and the function is not differentiable at the boundary point of this line segment.


\section{Appendix}

\noindent \textbf{Proof of Lemma~\ref{lem:concave}.} Using
(\ref{eq:alternative-esp-for-x}) and
(\ref{eq:distribution-of-first-jump}) we write
\begin{equation}\label{eq:first-term}
\begin{split}
\int_0^{t}\P\{s \leq
\sigma_1\}&k(\vec{x}(s,\pi))ds=\int_0^{t}e^{-\lambda_0
s}\cdot \E\left\{e^{-(\lambda_1-\lambda_0)(s-\Theta)^+}\right\}k(\vec{x}(s,\vec{\pi}))ds
\\&=c \int_0^{t}ds e^{-\lambda_0
s}\ E\{e^{-(\lambda_1-\lambda_0)(s-\Theta)^+}\}\frac{\E\{1_{\{s
\geq \Theta
\}}e^{-(\lambda_1-\lambda_0)(s-\Theta)}\}}{\E\{e^{-(\lambda_1-\lambda_0)(s-\Theta)^+}\}}
\\& = \frac{c \pi}{\lambda_1}(1-e^{-\lambda_1 t})+ c \sum_{i=1}^{n}\pi_i\int_0^{t} e^{-\lambda_0 s} \left(\int_0^{s}
e^{-(\lambda_1-\lambda_0)(s-u)}f_i(u) du \right)ds,
\end{split}
\end{equation}
where $f_i(\cdot)$ is the probability density function of
$\Theta$, given that $M_0=i$. Therefore,
$\vec{\pi} \rightarrow \int_0^{t}\P\{s \leq
\sigma_1\}\, k(\vec{x}(s,\vec{\pi}))$ is linear.
 Next, we observe that
\begin{equation}\label{eq:h-P}
\begin{split}
h(\vec{x}(t,\pi)) \, \P\{t<\sigma_1\}&=(1-x_{\Delta}(t,\vec{\pi}))\cdot e^{-\lambda_0
t} \cdot \E\left\{e^{-(\lambda_1-\lambda_0)(t-\Theta)^+}\right\}
\\&=e^{-\lambda_0 t} \cdot \E\left\{1_{\{t<\Theta\}}\right\}=e^{-\lambda_0
t}\sum_{i=1}^{n}\pi_i \int_t^{\infty}f_i(s)ds,
\end{split}
\end{equation}
hence the mapping
\begin{equation}\label{eq:second-term}
\vec{\pi}\rightarrow h(\vec{x}(t,\vec{\pi}))\P(t<\sigma_1) \quad
\text{is linear}.
\end{equation}
Finally, let $w(\cdot)$ be a positive and concave function. Then
it can be written as
\begin{equation}
w(\vec{\pi})=\inf_{k \in
K}\left(\beta_0^{(k)}+\beta_1^{(k)}\pi_1+\cdots+\beta_n^{(k)} \pi_n+\beta_{\Delta}^{(k)}\pi \right),
\end{equation}
for some index set $K$ and constants $\beta_j^{(k)}$. Using this
representation of $w(\cdot)$, (\ref{eq:alternative-esp-for-x}),
and (\ref{eq:sigma-density}) we obtain
\begin{equation} \label{eq:int-of-Sw-in-explicit-form}
\begin{split}
&\int_0^{t}\P(\sigma_1 \in ds)
Sw(\vec{x}(s,\vec{\pi})) \\ &=\int_0^{t}\P\{\sigma_1 \in
ds\}w\left(\frac{\lambda_0
x_1(s,\pi_1)}{\lambda_0(1-x_{\Delta}(s,\pi))+\lambda_1
x_{\Delta}(s,\pi)},\cdots , \frac{\lambda_1
x_{\Delta}(s,\pi_1)}{\lambda_0(1-x_{\Delta}(s,\pi))+\lambda_1
x_{\Delta}(s,\pi)}\right)
\\&=\int_0^{t}ds \, e^{-\lambda_0 s}\left[\lambda_0 \E\{1_{\{s<\Theta\}}\}+\lambda_1 \E\left\{1_{\{s \geq \Theta\}}
e^{-(\lambda_1-\lambda_0)(s-\Theta)}\right\}\right] \cdot
\\ &\qquad \qquad \qquad \inf_{k \in
K}\left[\beta_0^{(k)}+\frac{\beta_1^{(k)}\lambda_0
\E\{1_{\{M_s=1\}}\} +\cdots+\beta_{\Delta}\lambda_1\E\left\{1_{\{s
\geq
\theta\}}e^{-(\lambda_1-\lambda_0)(s-\Theta)}\right\}}{\lambda_0\E\left\{
1_{\{s<\Theta\}}\right\}+\lambda_1 \E\left\{1_{\{s \geq \Theta\}}
e^{-(\lambda_1-\lambda_0)(s-\Theta)}\right\} } \right]
\\&=\int_0^{t}ds \, e^{-\lambda_0 s} \bigg(\inf_{k \in
K} \beta_0^{(k)}\left(\lambda_0 \E\{1_{\{s<\Theta\}}\}+\lambda_1
\E\left\{1_{\{s \geq \Theta\}}
e^{-(\lambda_1-\lambda_0)(s-\Theta)}\right\}\right)
\\ &\qquad \qquad \qquad \qquad \qquad +\beta_1^{(k)}\lambda_0 \E\{1_{\{M_s=1\}}\}
+\cdots+\beta_{\Delta}\lambda_1\E\left\{1_{\{s \geq
\theta\}}e^{-(\lambda_1-\lambda_0)(s-\Theta)}\right\}\bigg).
\end{split}
\end{equation}
Note that the term inside the parentheses is linear in $\vec{\pi}$ due to (\ref{def:P-pi}). Hence it
 follows that
 \begin{equation}
\vec{\pi} \rightarrow \int_0^{t}\P\{\sigma_1 \in
ds\}Sw(\vec{x}(s,\vec{\pi})) \quad \text{is concave,}
 \end{equation}
since the lower envelope of linear functions is concave.

As a sum of three concave mappings, $\vec{\pi} \rightarrow
Jw(t,\vec{\pi})$ is concave for all $t \geq 0$. Also, as the lower
envelope of concave functions, the mapping $\vec{\pi} \rightarrow
J_0 w(\vec{\pi})=\inf_{t \geq 0}J (t,\vec{\pi})$ is again concave.
\hfill $\square$ \\

\noindent\textbf{Proof of Lemma~\ref{lem:cont}.}
Let $w(\cdot)$ be a bounded continuous function.
Then as in (\ref{eq:int-of-Sw-in-explicit-form}), we have
\begin{equation} \label{eq:int-of-Sw}
\begin{split}
\int_0^{t}\P(\sigma_1 \in ds) Sw(\vec{x}(s,\vec{\pi}))=\int_0^{t}
\lambda_0 e^{-\lambda_0 s} \left( \sum_{i=1}^{n} \pi_i
\int_s^{\infty}f_i(u)du \right) Sw(\vec{x}(s,\vec{\pi}))ds \, +
\\ \int_0^{t}\lambda_1 \left(\pi e^{-\lambda_1 s}+ e^{-\lambda_0 s} \sum_{i=1}^{n}\int_0^{s}e^{-(\lambda_1-\lambda_0)(s-u)}f_i(u) \right)
S w(\vec{x}(s,\vec{\pi}))ds.
\end{split}
\end{equation}
where $f_i(\cdot)$ is the probability density function of $\Theta$,
given that $M_0=i$. Then using (\ref{eq:first-term}), (\ref{eq:h-P}) and (\ref{eq:int-of-Sw}),
it can easily verified that the mapping $(t,\vec{\pi}) \rightarrow
Jw(t,\vec{\pi})$ is jointly continuous on $\R_+ \times D$. The mapping $(t,\vec{\pi}) \rightarrow
Jw(t,\vec{\pi})$ is then uniformly continuous on $[0,k]\times D$ for
all $k \in \mathbb{N}$. Therefore, the mapping
\begin{equation}\label{eq:J-k-cont}
\vec{\pi} \rightarrow J_{0,k}w (\vec{\pi}) = \inf_{t \in
[0,k]} J w(t,\vec{\pi})\quad  \text{is continuous on $D$}.
\end{equation}
On the other hand, using (\ref{eq:exp-for-J-w}) and (\ref{eq:h-P})
we can write
\begin{equation}\label{eq:bound-J-w}
\begin{split}
J w(t,\vec{\pi})&=Jw(t \wedge k, \vec{\pi})+\int_{t \wedge k
}^{t}\P\{s \leq \sigma_1\} k(\vec{x}(s,\pi))ds+
\\& e^{-\lambda_0 t}
\sum_{i\neq \Delta} \pi_i \left(\int_t^{\infty}f_i(s)ds- \int_{t
\wedge k }^{\infty}f_i(s)ds\right)+\int_{t \wedge
k}^{t}\P(\sigma_1 \in ds) Sw(\vec{x}(s,\vec{\pi}))
\\ &\geq Jw (t \wedge k,\vec{\pi})-e^{-\lambda_0 k} \sum_{i=1}^{n}\pi \int_{t \wedge k}^{\infty}f_i(s)ds
-  (\lambda_0 \vee \lambda_1) \int_{t \wedge k}^t  e^{- (\lambda_0 \wedge \lambda_1)s } \cdot ||w|| \, ds
\\ &\geq  Jw (t \wedge k,\vec{\pi})-e^{-(\lambda_0 \wedge \lambda_1) k}
\left(  1 + (\lambda_0 \vee \lambda_1) \cdot ||w||\right),
\end{split}
\end{equation}
where $||w|| \triangleq \sup_{\vp \in D} |w(\vp)|$. By taking the infimum on both sides of (\ref{eq:bound-J-w}) we get
\begin{equation}\label{eq:J-k-uniform}
J_{0,k}w (\vec{\pi}) \geq J_0 w(\vec{\pi}) \geq J_{0,k} w
(\vec{\pi}) -e^{-(\lambda_0 \wedge \lambda_1) k},
\end{equation}
which implies that $J_{0,k}w(\cdot) \rightarrow J_0 w(\cdot)$
uniformly on $D$. This fact together with
(\ref{eq:J-k-cont}) implies that $\vec{\pi} \rightarrow J_0 w
(\vec{\pi})$ is continuous on $D$. \hfill $\square$ \\

 \noindent \textbf{Proof of Proposition~\ref{prop:V-n-epsilon}.} First, we will prove (\ref{eq:eps-opt}) by an induction on $m \in
\mathbb{N}.$ For $m=1$ the left-hand-side of (\ref{eq:eps-opt})
becomes
\begin{equation} \label{prof-for-m-1}
\begin{split}
\E\left\{\int_0^{S_1^{\eps}}k(\vec{\Pi}_t)dt+h(\vec{\Pi}_{S_1^{\eps}})\right\}
&= \E \left\{\int_0^{r_0^{\eps}(\vec{\pi})\wedge
\sigma_1}k(\vec{\Pi}_t)dt+h(\vec{\Pi}_{r^{\eps}_0(\vec{\pi})\wedge
\sigma_1 })\right\}
\\ & =J v_0 (r_0^{\eps}(\vec{\pi}), \vec{\pi}) \leq J_0 v_0
(\vec{\pi})+\eps=v_1(\vec{\pi})+\eps,
\end{split}
\end{equation}
where we used (\ref{eq:defn-J}), (\ref{eq:seq-of-func}) and
(\ref{eq:defn-r-m}). Also note that we used Remark~\ref{rem:inf-J-is-attained} for the inequality above.
This inequality in (\ref{prof-for-m-1}) proves (\ref{eq:eps-opt}) holds for
$m=1$. Now, suppose (\ref{eq:eps-opt}) holds for $\eps \geq 0$,
and for some $m>1$. We will prove that it also holds when $m$ is
replaced by $m+1$. Since $S_{m+1}^{\eps} \wedge
\sigma_1=r_{m}^{\eps/2}(\vec{\pi}) \wedge \sigma_1$, we have
\begin{equation}\label{eq:derovation-f-m}
\begin{split}
&\E\left\{\int_{0}^{S_{m+1}^{\eps}}k(\vec{\Pi}_t)dt+
h(\vec{\Pi}_{S_{m+1}^{\eps}})\right\}
\\&=\E\bigg\{\int_{0}^{S_{m+1}^{\eps} \wedge \sigma_1
}k(\vec{\Pi}_t)dt+1_{\left\{S_{m+1}^{\eps} \geq \sigma_1
\right\}}\left[\int_{\sigma_1}^{S_{m+1}^{\eps}}k(\vec{\Pi}_t)dt+h(\vec{\Pi}_{S_{m+1}^{\eps}})\right]
+h(\vec{\Pi}_{S_{m+1}^{\eps}}) 1_{\{S_{m+1}^{\eps}<\sigma_1\}}
\bigg\}
\\&=\E\bigg\{\int_{0}^{r_{m}^{\eps/2}(\vec{\pi})
\wedge \sigma_1
}k(\vec{\Pi}_t)dt+1_{\left\{r_{m}^{\eps/2}(\vec{\pi}) \geq
\sigma_1 \right\}}\left[\int_{\sigma_1}^{\sigma_1+S_m^{\eps/2}
\circ \theta{\sigma_1}
}k(\vec{\Pi}_t)dt+h(\vec{\Pi}_{\sigma_1+S_{m}^{\eps/2}\circ
\theta_{\sigma_1} })\right]
\\&\qquad \qquad \qquad+h(\vec{\Pi}_{r_{m}^{\eps/2}(\vec{\pi}) \wedge \sigma_1})
1_{\{r_{m}^{\eps/2}(\vec{\pi})<\sigma_1\}} \bigg\}
\\&=\E\bigg\{\int_{0}^{r_{m}^{\eps/2}(\vec{\pi})
\wedge \sigma_1
}k(\vec{\Pi}_t)dt+h(\vec{\Pi}_{r_{m}^{\eps/2}(\vec{\pi}) \wedge
\sigma_1})
1_{\{r_{m}^{\eps/2}(\vec{\pi})<\sigma_1\}}\bigg\}+\E\left\{1_{\left\{r_{m}^{\eps/2}(\vec{\pi})
\geq \sigma_1 \right\}} f_{m}(\vec{\Pi}_{\sigma_1})\right\},
\end{split}
\end{equation}
in which
\begin{equation}
f_{m}(\vec{\pi})=\E\left\{\int_0^{S_n^{\eps/2}}k(\vec{\Pi}_t)dt+h(\vec{\Pi}_{S_n^{\eps/2}})\right\}
\leq v_{m}(\vec{\pi})+\eps/2,
\end{equation}
where the inequality follows from the induction hypothesis, and the
last line of (\ref{eq:derovation-f-m}) follows from the Strong
Markov property of the process $\vec{\Pi}$. Then we obtain
\begin{multline*}
\E\left\{\int_{0}^{S_{m+1}^{\eps}}k(\vec{\Pi}_t)dt+
h(\vec{\Pi}_{S_{m+1}^{\eps}})\right\}
\leq \E\bigg\{\int_{0}^{r_{m}^{\eps/2}(\vec{\pi})
\wedge \sigma_1
}k(\vec{\Pi}_t)dt+h(\vec{\Pi}_{r_{m}^{\eps/2}(\vec{\pi}) \wedge
\sigma_1})
1_{\{r_{n}^{\eps/2}(\vec{\Pi}_0)<\sigma_1\}}\bigg\}
\\ +\E\left\{1_{\left\{r_{n}^{\eps/2}(\vec{\Pi}_0)
\geq \sigma_1 \right\}}
v_{m}(\vec{\Pi}_{\sigma_1})\right\}+\frac{\eps}{2}
=J v_n (r_{m}^{\eps/2}(\vec{\pi}),\vec{\pi})+\frac{\eps}{2}
\leq v_{m+1}(\vec{\pi})+\eps,
\end{multline*}
where the first equality follows from the definition of the
operator $J$ in (\ref{eq:defn-J}) and the second equality follows
from (\ref{eq:defn-r-m}). This concludes the proof of
(\ref{eq:defn-S-eps}).

The inequality $V_n \leq v_n$ follows immediately from
(\ref{eq:eps-opt}) since $S_n^{\eps} \leq \sigma_n$ by
construction. Let us prove the opposite inequality $V_n \geq v_n$. First, we
will establish 
\begin{equation}\label{eq:V-vs-v}
\E\left\{\int_0^{\tau \wedge \sigma_m}k(\vec{\Pi}_t)dt+h
(\Vec{\Pi}_{\tau \wedge \sigma_m})\right\} \geq v_m(\vec{\pi}),
\end{equation}
for every $m \in \mathbb{N}$, by showing that
\begin{equation}\label{eq:V-vs-v-int}
\begin{split}
&\E\left\{\int_0^{\tau \wedge \sigma_m}k(\vec{\Pi}_t)dt+h
(\Vec{\Pi}_{\tau \wedge \sigma_m})\right\}
\\ & \geq \E\left\{\int_0^{\tau \wedge
\sigma_{m-k+1}}k(\vec{\Pi}_t)dt+1_{\{\tau \geq \sigma_{m-k+1}\}}
v_{k-1}(\vec{\Pi}_{\sigma_{m-k+1}})+1_{\{\tau <
\sigma_{m-k+1}\}}h(\vec{\Pi}_{\tau}) \right\}=:RHS_{k-1},
\end{split}
\end{equation}
for $k=1, \cdots, m+1$. Note that (\ref{eq:V-vs-v}) follows from
(\ref{eq:V-vs-v-int}) when we take $k=m+1$. For $k=1$, (\ref{eq:V-vs-v}) is satisfied as an equality since
$v_0(\cdot)=h(\cdot)$. Now, let us suppose that (\ref{eq:V-vs-v}) holds
for some $1 \leq k < m+1$, and let us prove that it also holds for $k+1$.

Note that $RHS_{k-1}$ can be decomposed as
\begin{equation}
RHS_{k-1}=RHS_{k-1}^{(1)}+RHS_{k-1}^{(2)},
\end{equation}
in which
\begin{equation}
\begin{split}
RHS_{k-1}^{(1)} & \triangleq \E\left\{\int_0^{\tau \wedge
\sigma_{m-k}}k(\vec{\Pi}_t)dt+
1_{\{\tau<\sigma_{m-k}\}}h(\vec{\Pi}_{\tau})\right\},
\\RHS_{k-1}^{(2)} & \triangleq \E\Bigg\{
1_{\{\tau \geq \sigma_{m-k}\}}\Big[\int_{\sigma_{m-k}}^{\tau
\wedge \sigma_{m-k+1}}k(\vec{\Pi}_t)dt
\\ &\qquad \qquad  +1_{\{\tau \geq
\sigma_{m-k+1} \}}v_{k-1}
(\vec{\Pi}_{\sigma_{m-k+1}})+1_{\{\tau<\sigma_{m-k+1}\}}h(\Pi_{\tau})\Big]\Bigg\}.
\end{split}
\end{equation}
By Lemma~\ref{lem:bremaud}, there exists an
$\F_{\sigma_{m-k}}$-measurable random variable $R_{m-k}$ such that
\begin{align*}
\tau \wedge \sigma_{m-k+1}=(\sigma_{m-k}+R_{m-k}) \wedge
\sigma_{m-k+1} \quad \text{ on } \; \{\tau \geq \sigma_{m-k}\}.
\end{align*}
Then
$RHS_{k-1}^{(2)}$ can be written as $RHS_{k-1}^{(2)}=$
\begin{equation}\label{eq:RHS-k-1-2}
\begin{split}
& \E\bigg\{ 1_{\{\tau \geq
\sigma_{m-k}\}}\bigg[\int_{\sigma_{m-k}}^{(\sigma_{m-k}+R_{m-k})
\wedge \sigma_{m-k+1}}k(\vec{\Pi}_t)dt+1_{\{\sigma_{m-k}+R_{m-k}
\geq \sigma_{m-k+1}\}}v_{k-1} (\vec{\Pi}_{\sigma_{m-k+1}})
\\ &+1_{\{\sigma_{m-k}+R_{m-k}<\sigma_{m-k+1}\}}h(\Pi_{\sigma_{m-k}+R_{m-k}})\bigg]\bigg\}
=\E\left\{1_{\{\tau \geq
\sigma_{m-k}\}}g_{m-k}(R_{m-k},\vec{\Pi}_{\sigma_{m-k}})\right\},
\end{split}
\end{equation}
in which
\begin{equation}
\begin{split}
g_{m-k}(r,\vec{\pi})& \triangleq \E\left\{\int_0^{r \wedge
\sigma_1}k(\vec{\Pi}_t)dt+1_{\{r \geq
\sigma_1\}}v_{k-1}(\Pi_{\sigma_1})+1_{\{r <
\sigma_1\}}h(\vec{\Pi}_{\tau}) \right\}
\\&=J v_{k-1}(r,
\vec{\pi}) \geq J_{0} v_{k-1}(\vec{\pi})=v_{k}(\vec{\pi}).
\end{split}
\end{equation}
The second equality in (\ref{eq:RHS-k-1-2}) follows from the
strong Markov property of the process $\Pi$ and the fact that the
jump times of the observation process $X$ and $\Pi$ are the same.
Therefore, the expression for $RHS_{k-1}^{(2)}$ is bounded below as
\begin{equation}
RHS_{k-1}^{(2)} \geq \E\left\{ 1_{\{\tau \geq \sigma_{m-k}\}}
v_{k}(\vec{\Pi}_{\sigma_{m-k}})\right\}.
\end{equation}
Therefore,
\begin{multline}
\E\Big\{\int_0^{\tau \wedge \sigma_m}k(\vec{\Pi}_t)dt + h
(\Vec{\Pi}_{\tau \wedge \sigma_m})\Big\} \\  \geq
\E\left\{\int_0^{\tau \wedge \sigma_{m-k}}
k(\vec{\Pi}_t)dt+1_{\{\tau < \sigma_{m-k}\}}h
(\Vec{\Pi}_{\tau})+1_{\{\tau \geq
\sigma_{m-k}\}}v_{k}(\vec{\Pi}_{m-k})\right\}
\end{multline}
This completes the proof of (\ref{eq:V-vs-v-int}) by induction.
Equation (\ref{eq:V-vs-v}) follows when we set $k=n+1$. Finally, taking the
infimum of both sides in (\ref{eq:V-vs-v}), we arrive at the
desired inequality $V_n \geq v_n$. \hfill $\square$ \\

\noindent \textbf{Proof of Proposition~\ref{prop:v-V}.} Using
Proposition~\ref{prop:sequential}, Corollary~\ref{cor:little-v-n}
and (\ref{prop:V-n-epsilon}) we obtain
\begin{equation}
v(\vp)=\lim_{m \rightarrow \infty}v_{m}(\vp)=\lim_{m \rightarrow
\infty}V_m(\vp)=V(\vp), \quad \vp \in D,
\end{equation}
which proves the first statement of the proposition. To prove the second statement, we note that
the sequence $\{v_{n}\}_{n \geq 1}$ is decreasing and
\begin{equation}\label{eq:V-J-0-V}
\begin{split}
&V(\vp)=v(\vp)=\inf_{m \geq 1}v_{m}(\vp)=\inf_{m \geq 1}\inf_{t
\geq 0}J v_{m-1}(t,\vp)=\inf_{t \geq 0}\inf_{m \geq 1}J
v_{m-1}(t,\vp)
\\ & =\inf_{t \geq 0}\inf_{m\geq 1} \Bigg\{\int_0^{t}\P\left\{s \leq \sigma_1\right\}k(\vx(s,\vp))ds
\\ &\qquad \qquad \qquad  \qquad \qquad +
\int_0^{t}\P\left\{\sigma_1 \in ds \right\}S
v_{m-1}(\vx(s,\vp))+\P\left\{t<\sigma_1\right\}h(\vx(t,\vp))\Bigg\}
\\&=\inf_{t \geq 0}\Bigg\{\int_0^{t}\P\left\{s \leq
\sigma_1\right\}k(\vx(s,\vp))ds+\int_0^{t}\P\left\{\sigma_1 \in ds
\right\}S
v(\vx(s,\vp))+\P\left\{t<\sigma_1\right\}h(\vx(t,\vp))\Bigg\}
\\&=\inf_{t \geq 0}Jv(t,\vp)=J_0 v (\pi).
\end{split}
\end{equation}
This proves that $V$ is a solution of $U=J_0 U$. The third line
of (\ref{eq:V-J-0-V}) follows from the bounded convergence
theorem. Next, let $U$ be a solution of $U=J_0 U$ such that $ U \leq h$. Then by Remark~\ref{rem:monotone} we have $U=J_0 U
\leq J_0 h =v_1$. Now, suppose $U \leq v_m$ for some $m \geq 0$,
then $U=J_0U \leq J_0 v_m=v_{m+1}$. By induction, we conclude that
$U \leq v_m$, for all $m \geq 1$ and therefore $U \leq \lim_{m
\rightarrow \infty}v_m=v=V$. \hfill $\square$ \\

\noindent \textbf{Proof of Lemma~\ref{lem:dyn-p-J-t-J}.} Let us
fix a constant $u \geq t$, and $\vp \in D$. Then
\begin{multline}
\label{eq:Ju}
J w (u,\vp)= \E\left\{\int_0^{u \wedge \sigma_1}k(\vP_s)ds+
1_{\{u<\sigma_1\}}h(\vP_u)+ 1_{\{u \geq
\sigma_1\}}w(\vP_{\sigma_1}) \right\}
\\ =  \E\left\{\int_0^{t \wedge \sigma_1}k(\vP_s)ds+
1_{\{u<\sigma_1\}}h(\vP_u)+ 1_{\{u \geq
\sigma_1\}}w(\vP_{\sigma_1})
\right\}  +  \E\left\{1_{\{\sigma_1>t\}}\int_t^{u \wedge
\sigma_1}k(\vP_s)ds\right\}.
\end{multline}
On the event $\{\sigma_1>t\}$,we have $u \wedge
\sigma_1=t+\left\{(u-t)\wedge (\sigma_1 \circ \theta_t)\right\}$.
Therefore, the strong Markov property implies
\begin{equation}
\begin{split}\label{eq:int-u-t-w-wedge-sigma1}
&\E\left\{1_{\{\sigma_1>t\}}\int_t^{u \wedge
\sigma_1}k(\vP_s)ds\right\}=\E\left\{1_{\{\sigma_1>t\}}\mathbb{E}^{\vP_t}\left\{\int_0^{(u-t)\wedge
\sigma_1}k(\vP_s)ds\right\}\right\}
\\ &=\E\left\{1_{\{\sigma_1>t\}}\left[J w (u-t, \vP_t)-\E\left\{1_{\{u-t<\sigma_1\}}h(\vP_{u-t})+1_{\{u-t \geq \sigma_1\}}w(\vP_{\sigma_1})\right\}\right]\right\}
\\&=\P\{\sigma_1>t\}J(u-t,\vx(t,\vp))-\E\left\{1_{\{u<\sigma_1\}}h(\vP_u)\right\}-\E\left\{1_{\{\sigma_1>t\}}
1_{\{u \geq \sigma_1\}}w(\vP_{\sigma_1})\right\},
\end{split}
\end{equation}
where the second equality follows from the definition of the
operator $J$, and the third from (\ref{eq:rel-pi-x}) and the
strong Markov property. Substituting
(\ref{eq:int-u-t-w-wedge-sigma1}) into (\ref{eq:Ju}), after some
simplification, yields
\begin{equation}
J w (u, \pi)=Jw (t,\pi)+\P\left\{\sigma_1>t\right\}\left[J
w(u-t,\vx(t,\vp))-h(\vx(t,\vp))\right].
\end{equation}
Now, taking the infimum of both sides over $u \in [t,\infty]$
concludes the proof. \hfill $\square$ \\

\noindent \textbf{Proof of Corollary~\ref{cor:r-m}.}
Note that by Remark~\ref{rem:inf-J-is-attained}, we have
\begin{equation}
Jv_{m}(r_m(\vp),\vp)=J_0v_{m}(\vp)=J_{r_{m}(\vp)}v_{m}(\vp).
\end{equation}
Let us first assume that $r_m(\vp)<\infty$. Taking $t=r_{m}(\vp)$
and $w=v_{m}$ in (\ref{eq:J-t-J}) gives
\begin{multline*}
Jv_{m}(r_m(\vp),\vp)=J_{r_{m}(\vp)}v_{m}(\vp) \\ =Jv_{m}(r_m(\vp),\vp)+\P\left\{\sigma_1>r_{m}(\vp)\right\}
\left[v_{m+1}(\vx(r_m(\vp), \vp))-h(\vx(r_m(\vp),\vp))\right].
\end{multline*}
Hence, we have $v_{m+1}(\vx(r_m(\vp),\vp))=h(\vx(r_m(\vp),\vp))$.

If $0<t<r_m(\vp)$, then
\begin{equation}\label{eq:Jvm-ge-J0}
Jv_m(t,\vp)>J_0 v_m(\vp)=J_{r_m(\vp)}v_m(\vp)=J_t v_m(\vp).
\end{equation}
Using (\ref{eq:J-t-J}) one more time, we get
\begin{equation}
J_0 v_m(\vp)=J_t v_m(\vp)=J v_m(t,\vp)+\P
\left\{\sigma_1>t\right\}\left[v_{m+1}(\vx(t,\vp))-h(\vx(t,\vp))\right].
\end{equation}
This equation together with (\ref{eq:Jvm-ge-J0}) implies that
$v_{m+1}(\vx(t,\vp))<h(\vx(t,\vp))$ for $t \in (0,r_m(\vp))$.

If $r_m(\vp)=\infty$, then $v_{m+1}(\vx(t,\vp))<h(\vx(t,\vp))$ for
every $t \in (0,\infty)$ by the same argument as in the last
paragraph. The statement of the lemma still holds in this case,
since by convention $\inf \emptyset=\infty$. \hfill $\square$ \\

\noindent \textbf{Proof of Proposition~\ref{prop:L-V}.} The proof
will be based on an induction. For $m=1$, by
Lemma~\ref{lem:bremaud} there exists a constant $u \in [0,\infty]$
such that $U_{\eps} \wedge \sigma_1= u\wedge \sigma_1$. Then
\begin{equation}\label{eq:E-L}
\begin{split}
&\E\left\{L_{U_{\eps}\wedge \sigma_1}\right\}=\E\left\{\int_0^{u
\wedge \sigma_1}k(\vP_s)ds+V(\vP_{u \wedge \sigma_1})\right\}
\\& =\E\left\{\int_0^{u \wedge \sigma_1}k(\vP_s)ds+1_{\{u \geq \sigma_1\}}V(\vP_{\sigma_1})+1_{\{u<\sigma_1\}}
h(\vP_u)\right\}+
\E\left\{1_{\{u<\sigma_1\}}[V(\vP_u)-h(\vP_u)]\right\}
\\&=J
V(u,\vp)+\P\left\{u<\sigma_1\right\}[V(\vx(u,\vp))-h(\vx(u,\vp))]=J_u
V(\vp),
\end{split}
\end{equation}
where the third equality follows from (\ref{eq:defn-J}) and
(\ref{eq:rel-pi-x}). The last equality follows from
(\ref{eq:dyn-prog-0.5}).

Fix any $t \in [0,u)$. By (\ref{eq:dyn-prog-0.5}) again
\begin{equation}
\begin{split}
&J
V(t,\vp)=J_tV(\vp)-\P\left\{\sigma_1>t\right\}[V(\vx(t,\vp))-h(\vx(t,\vp))]
\\ &\geq J_0
V(\vp)-\P\left\{\sigma_1>t\right\}[V(\vx(t,\vp))-h(\vx(t,\vp))]=J_0
V(\vp)-\E\left\{1_{\{\sigma_1>t\}}[V(\vP_t)-h(\vP_t)]\right\}.
\end{split}
\end{equation}
On the event $\{\sigma_1>t\}$ we have $U_{\eps}>t$ (otherwise,
$U_{\eps} \leq t <\sigma_1$ would imply $U_{\eps}=u \leq t$, and this 
would contradict our initial choice of $t<u$). Thus,
$V(\vP_t)-h(\vP_t)<-\eps$ on $\{\sigma_1>t\}$. Hence,
\begin{equation}
J V(t,\vp) \geq J_0 V(\vp)+\eps \P\left\{\sigma_1>t\right\}> J_{0}
V(\vp), \quad t \in [0,u).
\end{equation}
Therefore, $J_0 V(\vp)=J_u V(\pi)$ and (\ref{eq:E-L}) implies that
\begin{equation}
\E\left\{L_{U_{\eps}\wedge \sigma_1}\right\}=J_u V(\vp)=J_0 V
(\vp)=V(\vp)= L_0,
\end{equation}
which completes the proof for $m=1$.

Assume (\ref{eq:V-V}) holds for $m \geq 1$. Note that
\begin{equation}
\begin{split}
&\E\left\{L_{U_{\eps}\wedge
\sigma_{m+1}}\right\}=\E\left\{1_{\{U_{\eps}<\sigma_1\}}L_{U_{\eps}}+1_{\{U_{\eps}
\geq \sigma_1 \}}L_{U_{\eps}\wedge \sigma_{m+1}}\right\}
=\E\left\{1_{\{U_{\eps}<\sigma_1\}}L_{U_{\eps}}\right\}
\\&+\E\left\{1_{\{U_{\eps}
\geq \sigma_1 \}}\left[\int_{\sigma_1}^{U_{\eps}\wedge
\sigma_{m+1}}k(\vP_s)ds+V(\vP_{U_{\eps}\wedge
\sigma_{m+1}})\right]\right\}+\E\left\{1_{\{U_{\eps} \geq \sigma_1
\}} \int_0^{\sigma_1}k(\vP_s)ds\right\}.
\end{split}
\end{equation}
Since $U_{\eps}\wedge \sigma_{m+1}=\sigma_1+[(U_{\eps} \wedge
\sigma_m)\circ \theta_{\sigma_1}]$ on the event $\{U_{\eps} \geq
\sigma_1 \}$, the strong Markov property of $\vP$ implies that
\begin{multline}
\E\left\{L_{U_{\eps}\wedge
\sigma_{m+1}}\right\}=\E\left\{1_{\{U_{\eps}<\sigma_1\}}L_{U_{\eps}}\right\}
\\+
\E\left\{1_{\{U_{\eps} \geq \sigma_1
\}}\mathbb{E}^{\vP_{\sigma_1}}\left\{\int_0^{U_{\eps}\wedge
\sigma_m}k(\vP_s)ds+ V(\vP_{U_{\eps}\wedge
\sigma_m})\right\}\right\}
+\E\left\{1_{\{U_{\eps} \geq \sigma_1
\}} \int_0^{\sigma_1}k(\vP_s)ds\right\}.
\end{multline}
By induction hypothesis we can replace the inner expectation with
$V(\vP_{\sigma_1})$ and obtain
\begin{equation}
\begin{split}
\E\left\{L_{U_{\eps}\wedge
\sigma_{m+1}}\right\}&=\E\left\{1_{\{U_{\eps}<\sigma_1\}}L_{U_{\eps}}+1_{\{U_{\eps}
\geq \sigma_1
\}}\left[\int_{0}^{\sigma_1}k(\vP_s)ds+V(\vP_{\sigma_1})\right]\right\}
\\&=\E\left\{1_{\{U_{\eps}<\sigma_1\}}L_{U_{\eps}}+1_{\{U_{\eps}
\geq \sigma_1
\}} L_{\sigma_1} \right\}
=\E\left\{L_{U_{\eps}\wedge
\sigma_1}\right\}=L_0,
\end{split}
\end{equation}
where the last equality follows from the above proof for $m=1$.
This completes the proof of the statement. \hfill $\square$ \\

\bibliographystyle{dcu}
\bibliography{references}

\end{document}